\documentclass{article}

\usepackage[utf8]{inputenc}
\usepackage[OT2]{fontenc}
\usepackage[T1]{fontenc}
\usepackage{palatino}

\usepackage[english]{babel}
\usepackage{amsmath}
\usepackage{amsfonts}
\usepackage{mathrsfs}
\usepackage{amssymb}
\usepackage{stmaryrd}
\usepackage{chemfig}
\usepackage{mathtools}
\usepackage{amsthm}
\usepackage{hyperref}
\usepackage{wasysym}
\usepackage{enumerate}
\usepackage{xcolor}
\usepackage{float}

\usepackage{tikz-cd}
\usetikzlibrary{decorations.markings}
\usetikzlibrary{shapes,positioning,intersections,quotes}

\theoremstyle{plain}
\newtheorem{theo}{Theorem}[section]
\newtheorem{prop}[theo]{Proposition}
\newtheorem{coro}[theo]{Corollary}
\newtheorem{mylemma}[theo]{Lemma}

\theoremstyle{definition}
\newtheorem{madef}[theo]{Definition}

\newtheorem*{nota*}{Notation}
\newtheorem{rmq}[theo]{Remark}
\newtheorem{ex}[theo]{Example}
\newtheorem{discussion}{Discussion}

\theoremstyle{remark}
\newtheorem{myrmk}[theo]{Remark}

\definecolor{mypink2}{RGB}{219, 48, 122}

\newcommand{\red}[1]{\textcolor{black}{#1}}
\newcommand{\blue}[1]{\textcolor{black}{#1}}

\renewcommand{\P}{\mathcal{P}}
\renewcommand{\C}{\mathcal{C}}
\renewcommand{\O}{\mathbf{O}}
\newcommand{\Com}{\mathbf{Com}}
\newcommand{\ComP}{\Com^{\P}}
\newcommand{\Lie}{\mathbf{Lie}}
\newcommand{\Gerst}{\mathbf{Gerst}}
\newcommand{\LD}{\mathrm{LD}}
\newcommand{\tr}{\mathsf{TR}}

\newcommand{\CH}{\mathbf{CH}}
\newcommand{\GL}{\mathcal{GL}}
\newcommand{\GGL}{\Gerst}

\newcommand{\LBS}{\mathfrak{LBS}}
\renewcommand{\L}{\mathcal{L}}
\newcommand{\B}{\mathcal{B}}

\newcommand{\dg}{\mathrm{dg}}
\renewcommand{\d}{\mathrm{d}}
\renewcommand{\Bar}{\mathrm{\mathbf{B}}}
\newcommand{\F}{\mathcal{F}}
\newcommand{\Adm}{\mathrm{Adm}}
\newcommand{\E}{\mathcal{E}}
\newcommand{\Espec}{\mathrm{E}}

\newcommand{\FL}{\mathcal{FL}}
\renewcommand{\H}{\mathcal{H}}
\newcommand{\Kos}{\mathbf{Kos}}

\newcommand{\OS}{\mathrm{OS}}
\newcommand{\OSop}{\mathbf{OS}}

\newcommand{\zero}{\hat{0}}
\newcommand{\un}{\hat{1}}
\newcommand{\rk}{\mathrm{rk}\,}
\renewcommand{\L}{\mathcal{L}}
\newcommand{\EL}{\mathrm{EL}}

\newcommand{\K}{\mathbb{K}}

\newcommand{\Id}{\mathrm{Id}}

\newcommand{\op}{\mathrm{op}}

\newcommand{\lt}{\mathrm{lt}}
\newcommand{\G}{\mathcal{G}}

\newcommand{\N}{\mathbb{N}}
\newcommand{\Cx}{\mathbb{C}}
\newcommand{\Z}{\mathbb{Z}}

\newcommand{\kos}{\text{!`}}

\newcommand{\Q}{\mathbb{Q}}

\newcommand{\RKLS}{\mathbf{RKLS}}
\newcommand{\LKLS}{\mathbf{LKLS}}
\newcommand{\iRKLS}{\widehat{\mathbf{RKLS}}}
\newcommand{\iLKLS}{\widehat{\mathbf{LKLS}}}

\newcommand{\argmin}{\mathrm{argmin}\,}
\newcommand{\gr}{\mathrm{gr}\,}
\newcommand{\filtr}{\mathcal{F}}

\renewcommand{\d}{\mathrm{d}}

\newcommand{\iP}{P^{\circ}}
\newcommand{\nosP}{\P^{\star}}
\newcommand{\Hilb}{\mathrm{Hilb}}

\newcommand{\Eulerchar}{\mathrm{Euler}}

\begin{document}

\title{Operadic Kazhdan--Lusztig--Stanley theory}
\author{Basile Coron}
\date{}
\maketitle
\begin{abstract}
We introduce a new type of operad-like structure called a $\P$-operad, which depends on the choice of some collection of posets $\P$, and which is governed by chains in posets of $\P$. We introduce several examples of such structures which are related to classical poset theoretic notions such as poset homology, Cohen--Macaulayness and lexicographic shellability. We then show that $\P$-operads form a satisfactory framework to categorify Kazhdan--Lusztig polynomials of geometric lattices and their $P$-kernel. In particular, this leads to a new proof of the positivity of the coefficients of Kazhdan--Lusztig polynomials of geometric lattices. 
\end{abstract}

\section{Introduction}

The theory of Kazhdan--Lusztig--Stanley polynomials was introduced by \cite{Stanley_1992} in an attempt to unify the story of Kazhdan--Lusztig polynomials associated to Coxeter groups (\cite{KL_1979}) and the story of $g$-polynomials associated to polytopes (\cite{stanley_generalized_1987}) from a purely combinatorial standpoint. This framework would later be seen to encompass similar ``Kazhdan--Lusztig--like'' polynomials associated to other combinatorial objects such as matroids for instance (\cite{elias_kazhdanlusztig_2016}). The theory of Kazhdan--Lusztig--Stanley polynomials revolves around the key notion of a $P$-kernel.
\begin{madef}[P-kernel, \cite{Stanley_1992}]
Let $P$ be a locally finite graded bounded poset $P$. A $P$-\textit{kernel} $\kappa$ is a collection of polynomials $(\kappa_{XY})_{X\leq Y\in P}$ such that we have  
\begin{itemize}
    \item $\kappa_{XX}(t)= 1$ for all $X$ in $P$. 
    \item $\deg \kappa_{XY} \leq \rk [X,Y]$ for all $X\leq Y \in P$. 
    \item $\sum_{X\leq Y \leq Z}t^{\rk [X,Y]}\kappa_{XY}(t^{-1})\kappa_{YZ}(t) = 0 $ for all $X<Z \in P$. 
\end{itemize}
\end{madef}
From a $P$-kernel one can construct two polynomials $f$ and $g$, which we will call respectively left and right Kazhdan--Lusztig--Stanley polynomials (following \cite{brenti_twisted_1999}) via the following theorem. 
\begin{theo}[\cite{Stanley_1992}, Corollary 6.7]
Let $P$ be a locally finite graded bounded poset $P$ and $\kappa$ a $P$-kernel. There exists a unique collection of polynomials $(f_{XY})_{X\leq Y \in P}$ (resp. $(g_{XY})_{X\leq Y \in P}$) such that we have 
\begin{itemize}
    \item $f_{XX}(t) = 1$ (resp. $g_{XX}(t) = 1$) for all $X$ in $P$.
    \item $\deg f_{XY}(t) < \rk [X,Y] /2$ (resp. $\deg g_{XY}(t) < \rk [X,Y] /2$) for all $X<Y \in P$.
    \item $t^{\rk[X,Z]}f_{XZ}(t^{-1}) = \sum_{X\leq Y \leq Z}f_{XY}(t)\kappa_{YZ}(t)$ (resp. $t^{\rk[X,Z]}g_{XZ}(t^{-1}) = \sum_{X\leq Y \leq Z}\kappa_{XY}(t)g_{YZ}(t)$) for all $X<Z \in P$.  
\end{itemize}
\end{theo}
If $\kappa$ is a P-kernel then the collection $(t^{\rk [X,Y]}\kappa_{XY}(t^{-1}))_{XY}$ is also a $P$-kernel, whose KLS polynomials are called inverse left/right KLS polynomials associated to $\kappa$.

If $P$ is a Coxeter group with its Bruhat order, its $R$-polynomial (see  \cite{bjorner_combinatorics_2005} Chapter 5) is a $P$-kernel. The corresponding right KLS polynomial is the classical Kazhdan--Lusztig polynomial defined in \cite{KL_1979}. 

If $P$ is an Eulerian poset, then the collection of polynomials $\kappa$ defined by $\kappa_{XY}(t) \coloneqq (t - 1)^{\rk [X, Y]}$ for all $X \leq Y \in P$ is a $P$-kernel (\cite{Stanley_1992}, Proposition 7.1). In the case where $P$ is the face lattice of a polytope $\Delta$, the corresponding left KLS polynomial is the $g$-polynomial of $\Delta$.

For any locally finite graded bounded poset $P$, the characteristic polynomial of each interval of $P$ is a $P$-kernel (\cite{Stanley_1992}, Example 6.8). If $P$ is a geometric lattice then the corresponding right KLS polynomial is the Kazhdan--Lusztig polynomial of $P$ introduced in \cite{elias_kazhdanlusztig_2016}. All the KLS polynomials cited above were proved to have non-negative coefficients over the last two decades, via a common heuristical slogan: ``KLS polynomials are the Poincaré series of some stalk of some intersection cohomology sheaf''. In the realizable case (that is when the combinatorial object comes from a geometric object) this slogan is to be taken quite literally (see \cite{kazhdan_schubert_1980} for finite Weyl groups, \cite{stanley_generalized_1987} for rational polytopes and \cite{elias_kazhdanlusztig_2016} for arrangements). \cite{proudfoot_algebraic_2018} showed that in this case those three results can be unified under a common geometric framework. In the non-realizable case however, suddenly short of two millenia of geometry one has to rebuild a suitable cohomological theory from scratch, which is a daunting task (see \cite{Elias_2014} for Coxeter groups, \cite{karu_hard_2002} for non-rational polytopes, and \cite{braden2023singular} for geometric lattices). In this article we propose an alternative way of categorifying KLS polynomials of geometric lattices, which is not based on cohomological heuristics and instead relies on more involved global algebraic structures, which one could call ``operadic''. \\

Before considering the KLS polynomials of geometric lattices, let us focus on their corresponding $P$-kernel, i.e. the characteristic polynomials $\chi_P(t) \coloneqq \sum_{G \in P}\mu([\zero, G])t^{\rk [G, \un]}$. Those polynomials have no chance of being categorifiable as they can have negative coefficients. However, the polynomial $\chi_P^{+}(t) \coloneqq \sum_{G \in P}|\mu([\zero, G])|t^{\rk [G, \un]}$ (which differs from $\chi_P$ only by an alternating sign) is known to be the Poincaré series of a graded commutative algebra called the Orlik--Solomon algebra of $P$, denoted $\OS(P)$. If $P$ is realized by some complex hyperplane arrangement $\H$ the Orlik--Solomon algebra of $P$ is classically known to be isomorphic to the cohomology algebra over $\Q$ of the arrangement complement of $\H$ (see \cite{OS_1980}). For instance, if $P$ is the lattice of partitions of $\{1,\ldots, n\}$, then $P$ is realized by the braid arrangement $\{z_i = z_j, i,j \leq n\}$ and thus $\OS(P)$ is isomorphic to the cohomology algebra of the $n$-configuration space of $\Cx$. Alternatively, in this particular case $\OS(P)$ is also isomorphic to the cohomology algebra of the space of $n$-configurations of $2$-discs inside the unit $2$-disc (denoted $\LD_2(n)$), as this latter space is homotopically equivalent to the $n$-configuration space of $\Cx$. The collection of spaces $\{\LD_2(n), n \in \N^{\star}\}$ is known to have an interesting global associative structure called an operadic structure, which consists primarily of the maps
\begin{equation*}
    \LD(p)\times \LD(q) \xrightarrow{\circ_i} \LD(p+q-1)\, \,  \, \, i\leq p 
\end{equation*}
(the so-called ``operadic products'') given by inserting configurations of discs inside the $i$-th disc of a configuration of discs (an operation which was not possible with configurations of points). We refer to \cite{LV_2012} for a general reference on operads. This operadic structure is referred to as the little (2-)discs operad and is a cornerstone of operadic theory (see \cite{May1972} for more details on this central object). Those operadic products induce morphisms at the level of homology over $\Q$ (i.e. linear dual of Orlik--Solomon algebras) which form an operad in graded cocommutative coalgebras called $\Gerst$ (\cite{Cohen_1973}) which encodes Gerstenhaber algebras (\cite{gerst_alg_1963}). This operad satisfies a lot of nice properties, one of which being that it is Koszul (a property of operads parallel to the namesake property for associative algebras, see \cite{polishchuk_quadratic_2005}). This means that the Koszul complex of $\Gerst$ is acyclic. Finally, it turns out that the Euler characteristic of the Koszul complex of $\Gerst$ in arity $n$ is exactly the polynomial $\sum_{\zero \leq Y \leq \un \in \Pi_n }t^{\rk [\zero,Y]}\chi_{\zero Y}(t^{-1})\chi_{Y\un}(t)$ (the loss of sign due to considering $\chi^{+}$ instead of $\chi$ is compensated by the signs coming from the Euler characteristic) which recovers the fact that the characteristic polynomials of the intervals of a partition lattice form a kernel for that partition lattice. This hints at a connection between operadic theory (more specifically Koszulness of operads) on the one hand and KLS theory on the other. \\

In \cite{coron2023matroids} we proved that the operadic structure highlighted previously for partition lattices can be extended to the whole collection of geometric lattices, which gives a structure axiomatically similar to an operad (but much bigger), called an $\LBS$-operad. In this paper we generalize the notion of an $\LBS$-operad from the collection of geometric lattices to any collection $\P$ of finite bounded posets stable under taking closed intervals. We call the corresponding algebraic structure a $\P$-operad \blue{(Definition \ref{defpcoll}), which is roughly a collection of objects indexed by $\P$ together with morphisms satisfying an axiom resembling associativity.} In the first four sections of this paper we develop the ground theory of $\P$-operads, mostly showing that it is similar to that of associative algebras, with familiar notions such as presentations of $\P$-operads (Definition \ref{defpres}), Gröbner bases for $\P$-operads (Section \ref{secgrobner}) and finally Koszulness of $\P$-operads (Section \ref{seckos}) via either bar constructions (Definition \ref{defbarcons}) or Koszul complexes (Definition \ref{defkoscomplex}). A recurring toy example which we will use to illustrate those notions is given by the $\P$-operad in abelian groups $\Com^{\P}(P) \coloneqq \Z$ for all $P$ in $\P$, together with trivial operadic products. The notation comes from the analogy with the classical operad $\Com$ which encodes commutative algebras. For this particular example we show that the operadic notions cited above are closely connected to classical poset theoretic notions such as lexicographic shellability (Proposition \ref{propgrobcom}) and order complexes/poset homology (Remark \ref{exordbar}). \\

The last section of this article is devoted to applying the theory of $\P$-operads to the question of categorifying KLS polynomials of geometric lattices. For this application we will only use $\P$-operads with $\P$ the collection of geometric lattices. It is the author's hope that $\P$-operads for different collections $\P$ could be applied to other KLS polynomials. In Definition \ref{defKLScomplex} we introduce several graded differential complexes denoted $\RKLS$, $\LKLS$, $\iRKLS$ and $\iLKLS$ which are all subcomplexes of the bar construction of the operad of (linear duals of) Orlik-Solomon algebras. Those complexes are constructed so that their Euler characteristics in each grading  give the coefficients of \red{right} KLS polynomials, \red{left} KLS polynomials, inverse \red{right} KLS polynomials and inverse \red{left} KLS polynomials respectively (up to an alternating sign). We then prove the following theorem, which directly gives the claimed categorification. 

\begin{theo}[Theorem \ref{theomainhomo}]
Let $\L$ be a geometric lattice of rank $k$. 
\begin{enumerate}[i)]
    \item For $i<k/2$ the cohomology of $\RKLS_{(i)}(\L)$, $\LKLS_{(i)}(\L)$, $\iRKLS_{(i)}(\L)$ and $\iLKLS_{(i)}(\L)$ is concentrated in degree $i$. For $i> k/2$ the cohomology of $\RKLS_{(i)}(\L)$, $\LKLS_{(i)}(\L)$, $\iRKLS_{(i)}(\L)$ and $\iLKLS_{(i)}(\L)$ is concentrated in degree $i-1$. 
    \item If $k$ is even, the complexes $\RKLS_{(\frac{k}{2})}(\L)$, $\LKLS_{(\frac{k}{2})}(\L)$, $\iRKLS_{(\frac{k}{2})}(\L)$ and $\iLKLS_{(\frac{k}{2})}(\L)$ are acylic. 
\end{enumerate}
\end{theo}
We are mostly interested in the first part of Statement i), but one cannot prove this statement without the others. 
\section{$\P$-operads: first definitions and examples}\label{secdefex}
Throughout this section, let $\P$ denote a collection of finite bounded posets stable under taking closed intervals, fbs collection for short (e.g. all finite bounded posets, geometric lattices (see \cite{welsh_matroid_1976} for the definition), face lattices of polytopes, closed intervals of a given locally finite poset, etc). Let $\C$ be a monoidal category (e.g. modules over a ring with usual tensor product, $\dg$-modules over a ring with usual tensor product, commutative algebras with usual tensor product, topological spaces with cartesian product, etc, see \cite{mac_lane_categories_1978} Section VII for a general reference on monoidal categories). In this section we will introduce the notion of a $\P$-operad in $\C$ and develop the basic theory around it. At ground level this theory is very similar to that of associative algebras (\cite{bourbaki_algebra_1998}), classical operads (\cite{LV_2012}), and many other associative structures. We denote by $\nosP$ the collection of posets in $\P$ which are not singletons. For any bounded poset $P$ we denote by $\iP$ the subset $P\setminus\{\zero, \un\}$, called the interior of $P$.
\begin{madef}[$\P$-collection, $\P$-operad]\label{defpcoll}

    A $\P$-\textit{collection} $V$ in $\C$ is a collection $\{V(P), P \in \nosP\}$ of objects in $\C$ indexed by $\nosP$. A $\P$-\textit{operad} $\O = (O, \mu)$ consists of a $\P$-collection $O$ together with morphisms 
    \begin{equation*}
        \mu_{G,P} : O([\zero, G])\otimes O([G, \un]) \rightarrow O(P)
    \end{equation*}
    for any element $G \in \iP \subset P \in \P$, such that for any pair $G_1 < G_2$ of elements in $\iP \subset P \in \P$ we have the equality
    \begin{equation}\label{eqaxop}
        \mu_{G_2, P}\circ (\mu_{G_1, [\zero, G_2]} \otimes \Id) = \mu_{G_1, P}\circ (\Id \otimes \mu_{G_2, [G_1, \un]}).
    \end{equation}
    A $\P$-cooperad in $\C$ is a $\P$-operad in the opposite category $\C^{\op}$. In that case the morphisms will be denoted by $\Delta_{G,P}$ instead of $\mu_{G,P}$. 
\end{madef}
If the poset $P$ in which we are performing our operadic product $\mu_{G,P}$ is clear from the context we will omit it. Whenever we have a chain of elements $G_1 < \cdots < G_n $ in the interior of some poset $P \in \P$ one can compose operadic products $\mu_{G_i}$ in whichever order we prefer, to get a morphism: 
\begin{equation*}
    O([\zero, G_1])\otimes \cdots \otimes O([G_n , \un]) \rightarrow O(P).
\end{equation*}
Equality \eqref{eqaxop} implies  that this morphism does not depend on the order of composition we choose. We will denote this morphism $\mu_{G_1, \ldots, G_n, P}$. 
\begin{ex}\label{exPop}
\begin{itemize} 
    \item If $\P$ is the collection of products of partition lattices, a $\P$-operad is a classical object called a shuffle operad with levels (see \cite{DK_2010} for shuffle operads and \cite{Fresse_2003} for operads with levels), which is similar to a classical operad.  
    \item For any fbs collection $\P$ one can define the $\P$-operad $\Com^{\P}$ in abelian groups by $\mathrm{Com}(P) \coloneqq \Z$ for all $P \in \nosP$, and the obvious operadic products. The notation comes from the fact that when $\P$ is the collection of all products of partition lattices we get the shuffle operad with levels encoding commutative algebras.
    \item For any fbs collection $\P$ one can define a $\P$-collection $C_{\bullet}(-)$ in differential graded abelian groups by setting $C_{\bullet}(P)$ to be the simplicial chain complex associated to the order complex of $\iP$. This $\P$-collection has a cooperadic structure defined by 
    \begin{equation*}
    \begin{array}{ccl}
         C_{\bullet}(P) & \longrightarrow & \hspace{30px}C_{\bullet}([\zero, G])\otimes C_{\bullet}([G,\un])  \\
         \{ G_1 < \cdots < G_n \}  & \rightarrow & \left\{ \begin{array}{cc} \{G_1, \ldots, G_{i-1} \}\otimes \{G_{i+1},\ldots, G_n\} & \textrm{ if } G = G_i, \\
         0 & \textrm{ otherwise.} \end{array} \right. 
    \end{array}
    \end{equation*}
    \red{We shall see in Example \ref{exordbar} how this cooperad is related to the previous example.}
    \item Let $\GL$ denote the collection of geometric lattices. Any $\LBS$-operad, as introduced in \cite{coron2023matroids}, gives a $\GL$-operad when restricted to maximal building sets. This includes for instance the cooperad of Chow rings of matroids $\CH$ in the monoidal category of graded commutative algebras (the underlying collection being the Chow rings and the morphisms being induced by inclusions of torus orbit closures in the toric variety associated to the corresponding Bergman fan, see \cite{braden_semi-small_2022} for more details). Another example which will be of central importance in this article is given by the $\GL$-collection of Orlik--Solomon algebras (see \cite{Yuzvinsky_2001} for the definition). This collection admits \red{a co}operadic structure given by the morphisms of algebras $\Delta_G : \OS(\L) \rightarrow  \OS([\zero, G])\otimes \OS([G, \un]) $ defined \red{on the generators $e_H$ indexed by rank $1$ flats}, by
\begin{equation*}
\Delta_G(e_H) \coloneqq \left \{
\begin{array}{cccc}
    e_H&\otimes& 1 & \textrm{ if } H \leq G,  \\
    1 &\otimes & e_{G\vee H} & \textrm{ otherwise.}
\end{array}
\right.
\end{equation*}
Note that the very definition of this cooperadic structure uses geometricity in a very severe way. We refer to \cite{coron2023matroids} for more details. We will denote by $\OSop$ this cooperad. 

\end{itemize}
\end{ex}
\blue{An operad $\O = (O, \mu)$ in $\C$ will be called \textit{graded} if $\C$ is a monoidal category of graded objects (e.g. graded abelian groups with graded tensor product). Explicitly this means that each object $O(P)$ is graded and we have $$\mu_{G,P}(O_p([\zero, G])\otimes O_q([G, \un]))\subset O_{p+q}(P)$$ for all $P \in \nosP$ and all $G\in \iP$.} \blue{If $\C$ admits a duality functor $\vee$ (e.g. abelian groups or vector spaces over a given field) then one can define the dual $V^{\vee}$ of a $\P$-collection $V$ as the collection obtained by applying $\vee$ to each object of $V$.} If $\vee$ is strictly compatible with the monoidal product \red{(e.g. finitely generated abelian groups or finite dimensional vector spaces over a given field)} then the dual of a collection $O$ with operadic structure $\mu$ has a cooperadic structure given by $\mu^{\vee}$, and vice versa if we start with a cooperad. 
\begin{madef}[Morphism between $\P$-collections/operads]

Let $V,V'$ be two $\P$-collections.  A \textit{morphism} $\phi:V \rightarrow V'$ between $V$ and $V'$ is a collection of morphisms $\{\phi_P: V(P) \rightarrow V'(P), P \in \nosP \}$. 

Let $\O = (O, \mu)$ and $\O' = (O, \mu')$ be two $\P$-operads. A morphism $\phi:\O \rightarrow \O'$ between $\O$ and $\O'$ is a morphism between the underlying $\P$-collections $O$ and $O'$, which satisfies the compatibility relation
\begin{equation*}
    \phi_P\circ \mu_G = \mu'_G \circ (\phi_{[\zero, G]}\otimes \phi_{[G, \un]})
\end{equation*}
\blue{for all $P \in \nosP$ and all $G \in \iP.$}
\end{madef}
In the remainder of this section we assume that $\C$ is the monoidal category of $A$-modules over some ring $A$, with usual tensor product. 
\begin{madef}[Free operad generated by a collection]

Let $V$ be a $\P$-collection. The \textit{free} $\P$-\textit{operad} $\P(V)$ generated by $V$ is the operad consisting of the $\P$-collection
\begin{equation*}
    \P(V)(P) \coloneqq \bigoplus_{\zero < G_1 < \cdots < G_n <\un \subset P} V([\zero, G_1])\otimes \cdots \otimes V([G_n, \un])
\end{equation*}
with operadic morphisms $\mu_{G, P}: \P(V)([\zero, G])\otimes \P(V)([G,\un]) \rightarrow \P(V)(P)$ sending the summand 
\begin{multline*}
(V([\zero, G_1])\otimes \cdots \otimes V([G_n, G]))\otimes(V([G, G'_1])\otimes \cdots \otimes V([G'_{\blue{m}}, \un])) \subset \P(V)([\zero, G])\otimes \P(V)([G,\un])
\end{multline*} to the summand in $\P(V)(P)$ corresponding to the chain $\zero < G_1 <\cdots < G < G'_1 < \cdots < \un$ in $P$, via the identity. \red{The free operad $\P(V)$ is graded by the length of the chains. We will denote by $\P_k(V)$ the part of $\P(V)$ of grading $k$.}  
\end{madef}
We have an obvious inclusion of $\P$-collections $\iota_V: V \hookrightarrow \P(V)$. The terminology ``free'' is justified by the following straightforward proposition. 
\begin{prop}
Let $V$ be a $\P$-collection and $\O = (O, \mu)$ a $\P$-operad. For any morphism of $\P$-collections $\phi: V \rightarrow O$ there exists a unique morphism of $\P$-operads $\Tilde{\phi}: \P(V) \rightarrow \O$ such that we have $\Tilde{\phi} \circ \iota_V = \phi$. 
\end{prop}
In categorical terms $\P(-)$ is left-adjoint to the forgetful functor from $\P$-operads to $\P$-collections. 
\begin{madef}[Operadic ideal, ideal generated by a subcollection]

Let $\O = (O, \mu)$ be a $\P$-operad and $V = \{V(P) \subset O(P), P \in \P\}$ a subcollection of $O$. We say that $V$ is an \textit{ideal} of $\O$ if for any $P \in \P$ and any $G \in \iP$, the operadic product $\mu_{G, P}$ sends both $V([\zero, G])\otimes O([G, \un])$ and $O([\zero, G])\otimes V([G, \un])$ to $V(P)$. \\
We denote by $\langle V \rangle$ the subcollection of $O$ defined by 
\begin{equation*}
      \red{\langle V \rangle (P) \coloneqq \sum_{\zero < F < G < \un}\mu_{F,G,P}(O([\zero, F])\otimes V([F,G])\otimes O([G, \un])).}
\end{equation*}
One can check that $\langle V \rangle$ is the smallest ideal containing $V$ and we call it the ideal generated by $V$.
\end{madef}
\begin{ex}

For any morphism of operads $\phi: \O_1 \rightarrow \O_2$, the subcollection $\ker \phi \coloneqq \{\ker \phi_P \subset O_1(P), P \in \nosP\}$ is an ideal. 
\end{ex}

\begin{madef}[Operadic quotient]

Let $\O = (O, \mu)$ be a $\P$-operad and $I$ an ideal of $\O$. The \textit{operadic quotient} $\O/I$ is the $\P$-operad which consists of the $\P$-collection $$O/I(P) \coloneqq O(P)/I(P), \, \, \, \forall P \in \nosP$$ together with the operadic products \blue{
\begin{equation*}
    \frac{O([\zero,G])}{I([\zero,G])}\otimes \frac{O([G, \un])}{I([G, \un])} \simeq \frac{O([\zero, G])\otimes O([G, \un])}{I([\zero,G])\otimes O([G, \un]) \oplus O([\zero, G]) \otimes I([G, \un])} \longrightarrow \frac{O(P)}{I(P)}. 
\end{equation*}
induced by $\mu_{G,P}$ for all $P \in \nosP$ and all $G \in \iP$. }
\end{madef}
One can check that operadic quotients satisfy the usual universal property of quotients. 
\begin{madef}[Presentation of an operad]\label{defpres}

Let $\O= (O, \mu)$ be a $\P$-operad. A \textit{presentation} of $\O$ is the datum of a subcollection $V \subset O$ (the generators) such that the morphism $\P(V) \rightarrow \O$ \blue{induced by the inclusion $V\subset O$} is surjective over each poset $P \in \P$, together with a subcollection $R \subset \P(V)$ (the relations) such that $\langle R \rangle$ is the kernel of  the morphism $\P(V) \rightarrow \O$.
\end{madef}
This means that we have an isomorphism of $\P$-operads $\P(V)/\langle R \rangle \xrightarrow{\sim} \O.$ We say that a presentation $(V, R)$ is \blue{\textit{quadratic}} if $R$ is included in the part of $\P(V)$ of grading $2$. We say that an operad is \blue{\textit{quadratic}} if it admits a quadratic presentation. 
\begin{ex}

\begin{itemize}
\item Let us try to find a presentation of $\Com^{\P}$. One can see that the subcollection
\begin{equation*}
V(P) \coloneqq \left\{ 
\begin{array}{cl}
\Z & \textrm{ if } P \textrm{ has rank } 1, \\
\{0\} & \textrm{ otherwise,}
\end{array}
\right.
\end{equation*}
generates $\Com^{\P}$. Indeed for every poset $P$ in $\P$, by finiteness one can find a maximal chain $\zero = G_0 < G_1 < \cdots < G_n < G_{n+1} = \un$ in $P$, i.e. a chain such that every interval $[G_i, G_{i+1}]$ has rank $1$. This means that we have $\ComP([G_i, G_{i+1}]) = V([G_i, G_{i+1}])$ for all $i$. By definition of the operadic product in $\ComP$, the morphism $\mu_{G_1, \ldots , G_n}$ will send $1\otimes \cdots \otimes 1$ to $1 \in \ComP(P)$. In other words, the induced morphism $\P(V)\rightarrow \ComP$ is surjective. The question of the relations between those generators is more delicate and heavily depends on $\P$. One can readily see that the kernel of $\P(V) \rightarrow \ComP$ is linearly generated by the relations 
\begin{equation*}
    \mu_{G_1, \ldots, G_n}(1\otimes \cdots \otimes 1) - \mu_{G'_1, \ldots, G'_{m}}(1\otimes \cdots \otimes 1) \in \P(V)(P) 
\end{equation*}
where $G_1 < \cdots < G_n$ and $G'_1 < \cdots < G'_{m}$ run over maximal chains in $\iP \subset P \in \P$. However, what we are interested in is finding a smaller set of relations which operadically generates all the other relations. Ideally, we would like to limit ourselves to quadratic relations (i.e. relations between products over chains in posets of rank $2)$. For general $\P$, this is not possible but it is possible for many interesting collections of posets. For instance, let us show by hand that this is true when $\P$ is the collection $\GL$ of geometric lattices. Let $\L$ be a geometric lattice and $\zero = G_0 < \cdots < G_{n+1} = \un$, $\zero = G'_0 < \cdots < G'_{n+1} = \un$ be two maximal chains in $\L$ (since $\L$ is geometric those two chains must have the same cardinality). We want to show that the equality between operadic products over those two chains is a consequence of equalities between operadic products over chains in rank $2$. We will do that by induction on the rank of $\L$. If $G_1$ is equal to $G'_1$ then our induction hypothesis on $[G_1, \un]$ immediately gives us the result. Otherwise, let $G''_2$ be the join of $G_1$ \red{and} $G'_1$. By geometricity the element $G''_2$ has rank $2$. Let $G''_2 < G''_3 < \cdots < G''_{n+1} = \un$ be any maximal chain between $G''_2$ and $\un$. We have the equality
\begin{multline*}
    \mu_{G_1, \ldots, G_n}(1\otimes \cdots \otimes 1) - \mu_{G'_1, \ldots, G'_{n}}(1\otimes \cdots \otimes 1) =  \\ \mu_{G_1, \ldots, G_n}(1\otimes \cdots \otimes 1) -  \mu_{G_1, G''_2, \ldots, G''_n}(1\otimes \cdots \otimes 1) \\ + \mu_{G_1, G''_2, \ldots, G''_n}(1\otimes \cdots\otimes 1) - \mu_{G'_1, G''_2, \ldots, G''_n}(1\otimes\cdots \otimes 1) \\ + \mu_{G'_1, G''_2, \ldots, G''_n}(1\otimes \cdots \otimes 1) - \mu_{G'_1, \ldots, G'_{n}}(1\otimes \cdots \otimes 1). 
\end{multline*}
The first difference in the right hand term is quadratically generated, by our induction hypothesis on $[G_1, \un]$, the second difference is the operadic product of a quadratic relation (in $[\zero, G''_2]$) with the element $\mu_{G''_2, \ldots, G''_n}(1\otimes \cdots\otimes 1)$ and the last difference is quadratically generated, by our induction hypothesis on $[G'_1, \un]$. This concludes the proof. \\

As another \red{hands-on} example, let us consider the collection $\FL$ of face lattices of polytopes (sets of faces of a polytope ordered by inclusion), and let us show that $\Com^{\FL}$ is quadratic as well. As above consider two maximal chains $\zero < G_1 <\cdots < G_n <\un$ and $\zero <G'_1 < \cdots < G'_n < \un$. We cannot use the same trick as for geometric lattices because this time the join $G_1\vee G'_1$ may have rank strictly greater than $2$. However, since a polytope is connected, its $1$-skeleton must be connected as well, which means that there exists a sequence $G_1 = I_0 < I_1 > I_2 < \cdots > I_k = G'_1$ of covering relations in $P$. For each $I_p$ of rank $2$ choose any maximal chain from $I_p$ to $\un$ and repeat the same trick as for geometric lattices. In Section \ref{seckos} we will see that the common feature of $\FL$ and $\GL$ behind the quadraticity of $\Com$ is Cohen--Macaulayness (see \cite{bjorner_introduction_1982} for the definition). 
\item One can also define operads directly by giving a presentation. Here is an example, which will be of interest later on. For any $\P$ let us define the $\P$-collection $V$ by 
\begin{equation*}
V(P) \coloneqq \left\{ 
\begin{array}{cl}
\Z & \textrm{ if } P \textrm{ has rank } 1, \\
\{0\} & \textrm{ otherwise,}
\end{array}
\right.
\end{equation*}
and $R$ the subcollection of $\P(V)$ by 
\begin{equation*}
    R(P) \coloneqq \left \{
    \begin{array}{cl}
        \Z\langle \sum_{H\in \iP}\mu_{H, P}(1,1)\rangle & \textrm{ if } P \textrm{ has rank } 2,\\
    \{0\} & \mathrm{otherwise.} 
    \end{array}
    \right. 
\end{equation*}
We define the $\P$-operad $\Lie^{\P}$ in abelian groups by $ \Lie^{\P} \coloneqq \P(V)/\langle R \rangle.$ The notation comes from the fact that when $\P$ is the collection of all products of partition lattices, $\Lie^{\P}$ is the shuffle operad with levels encoding Lie algebras \blue{(see \cite{GK_1994})}. The unique relation in $R(\Pi_3)$ encodes the Jacobi identity. 
\end{itemize}
\end{ex}
In classical operadic theory, $\Com$ and $\Lie$ are related to each other by a notion called Koszul duality. In Section \ref{seckos} we will develop a theory of Koszulness for $\P$-operads and show that the Koszul duality between $\Com$ and $\Lie$ extends to this new context (see Example \ref{exdualliecom}). We are finally ready to define the main protagonist of Section \ref{secKLS}.

\begin{madef}[Gerstenhaber $\P$-operad]\label{defgerst}

For any $\P$ the $\P$-operad $\Gerst^{\P}$ in rational vector spaces is defined as the quotient $\P(V)/\langle R \rangle$ where $V$ is the $\P$-collection
\begin{equation*}
V(P) \coloneqq \left\{ 
\begin{array}{cl}
\Q C \oplus \Q L & \textrm{ if } P \textrm{ has rank } 1, \\
\{0\} & \textrm{ otherwise,}
\end{array}
\right.
\end{equation*}
with $C$ and $L$ two symbols (standing for Com and Lie respectively), and $R$ is the subcollection of $\P(V)$ defined by

\begin{align*}
R(P) \coloneqq \Q  &  \left\langle \mu_H(C\otimes L) - \sum_{H' \neq H}\mu_{H'}(L\otimes C),\,\,\,  H \in \iP, \right. \\
& \hspace{20px}\mu_{H}(C\otimes C) - \mu_{H'}(C\otimes C),\,\,\, H, H' \in \iP, \\ 
& \left. \hspace{40px}\sum_{H \in \iP} \mu_H(L\otimes L) \right\rangle
\end{align*}
if $P$ has rank $2$, and $\{0\}$ otherwise. \red{The defining relations of $\Gerst$ being homogeneous in both generators $C$ and $L$, each vector space $\Gerst^{\P}(P)$ \,($P \in \P$) is bigraded. For all integers $p,q$ we will denote by $\Gerst^\P_{(p,q)}$ the subcollection of $\Gerst^\P$ with $p$ the grading given by $C$ an $q$ the grading given by $L$.}
\end{madef}
The notation comes from the fact when $\P$ is the collection of products of partition lattices, $\Gerst^{\P}$ is the shuffle operad with levels encoding Gerstenhaber algebras (algebras with two binary products ($C$ and $L$), one of which is a Lie bracket (last relation of $R$), the other is a commutative product (second set of relations of $R$) and the Lie bracket is a derivation of the commutative product (first set of relations)). We refer to \cite{gerst_alg_1963} for more details and historical motivation behind Gerstenhaber algebras. \blue{In general, we shall see in Proposition \ref{propgrobgerst} that for every geometric lattice $\L$, the linar dual $\Gerst^{\GL}(\L)^{\vee}$ can be identified with the Orlik--Solomon algebra of $\L$.} \red{We could have defined $\Gerst$ over $\Z$ instead of $\Q$ (as for $\Lie$ and $\Com$) but in Section \ref{secKLS} we will need $\Gerst$ to be defined over a field of characteristic zero.}
\blue{\begin{ex}
If $P$ is of rank $1$ we obviously have $\Gerst^{\P}_{(1,0)}(P) = \Q C$ and $\Gerst^{\P}_{(0,1)}(P) = \Q L$. If $P$ is of rank $2$, we have 
\begin{equation*}
    \Gerst^{\P}_{(0,2)}(P) = \Q\langle \mu_H(L\otimes L), H \textrm{ atom of } P \rangle / \Q\sum_{H}\mu_H(L \otimes L), 
\end{equation*}
which has dimension the number of atoms of $P$ minus 1, or alternatively $\mu(P).$ We also have 
\begin{equation*}
    \Gerst^{\P}_{(2,0)}(P) = \Q\langle \mu_H(C\otimes C), H \textrm{ atom of } P \rangle / \langle \mu_H(C\otimes C) - \mu_{H'}(C \otimes C), H, H' \textrm{ atoms of } P \rangle  
\end{equation*}
which has dimension $1$. Finally, we have 
\begin{equation*}
    \Gerst^{\P}_{(1,1)}(P) = \frac{\Q\langle \mu_H(C\otimes L),  H \textrm{ atom of } P \rangle \oplus \Q\langle \mu_H(L\otimes C),  H \textrm{ atom of } P \rangle}{\langle \mu_H(C\otimes L) - \sum_{H' \neq H}\mu_{H'}(L\otimes C), \, H \textrm{ atom of } P\rangle }
\end{equation*}
which has dimension the number of atoms of $P$.
\end{ex}}
\section{Gröbner bases for $\P$-operads}\label{secgrobner}
Gröbner bases for algebras are a computational tool which is designed to deal with algebras defined by generators and relations. We refer to \cite{BW_1993} for a textbook reference on that subject. The general idea is to start by choosing a linearly ordered basis of the space of generators. This order is then used to derive an order on all monomials (pure tensors of elements of the basis in the free algebra), which is compatible in some sense with the multiplication of monomials (we call such orders ``admissible''). We then use this order to rewrite monomials in the quotient algebra:
\begin{equation*}
\textrm{greatest term } \longrightarrow \sum \textrm{lower terms},
\end{equation*}
for every relation $R = \textrm{greatest term} - \sum \textrm{lower terms}$ in some subset $\B$ of the ideal of relations (usually the greatest term is called the ``leading term'' and we will use this denomination). The subset $\B$ is called a Gröbner basis when it contains ``enough'' elements. To be precise we want that every leading term of some relation in the quotient ideal is divisible by the leading term of some element of $\B$. The general goal is to find a Gröbner basis as \red{small} as possible so that the rewriting is as easy as possible. At the end of the rewriting process (which stops if the monomials are well-ordered) we are left with all the monomials which are not rewritable i.e. which are not divisible by a leading term of some element of $\B$. Those monomials are called ``normal'' and they form a linear basis of our algebra exactly when $\B$ is a Gröbner basis. This basis comes with multiplication tables given by the rewriting process. \\

It turns out that this general strategy can be applied to structures which are much more general and complex than algebras, such as shuffle operads for instance (see \cite{DK_2010}). Loosely speaking, all we need in order to implement this strategy is to be able to make reasonable sense of the key words used above, such as ``monomials'', ``admissible orders'' and ``divisibility between monomials''. This is what we will set out to do for $\P$-operads in the following subsection. 

\subsection{First definitions}
Let $\P$ be an fbs collection and $V$ a $\P$-collection in the monoidal category of $A$-modules over some commutative ring $A$. Assume that we have chosen an additive basis $B(P)$ of $V(P)$ for each $P$ in $\P$. 

\begin{madef}[Monomial]
A \textit{monomial} in $\P(V)$ is an element of the form $\mu_{G_1, \ldots, G_n}(e_1\otimes \cdots \otimes e_{n+1}),$ for some chain $\zero < G_1 < \cdots < G_n < G_{n+1} = \un$ in some poset $P$ of $\P$, and some elements $e_i$  in $B([G_{i-1}, G_i])$ respectively. 
\end{madef}

In other words a monomial is given by a chain in some poset of $\P$, with each interval decorated by an element of the basis of generators over that interval. 

\begin{madef}[Divisibility between monomials]

A monomial $\mu_{G_1, \ldots, G_n}(e_1\otimes \cdots \otimes e_{n+1})$ in $\P(V)(P)$ is said to \textit{divide} another monomial $\mu_{G'_1, \ldots, G'_{n'}}(e'_1\otimes \cdots \otimes e'_{n'+1})$ if there exists $p$ and $q$ such that $P = [G'_p, G'_q]$, the two chains $G'_p < \cdots < G'_q$ and $ \blue{G'_p} < G_1 < \cdots < G_n < \blue{G'_q}$ are equal, and $(e_1, \ldots, e_{n+1}) = (e'_{p+1}, \ldots, e'_q)$. We extend this definition to terms of the form $\lambda m$ with $\lambda$ an element of the ring $A$ and $m$ a monomial, by saying that $\lambda_1 m_1$ divides $\lambda_2 m_2$ if $\lambda_1$ divides $\lambda_2$ and $m_1$ divides $m_2$.
\end{madef}

In plain English, a monomial divides another monomial if the chain of the dividing monomial is a subchain of the second and the decorations over this subchain coincide. 

\begin{madef}[Admissible ordering]

A \blue{total} ordering $\trianglelefteq$ on monomials is said to be \textit{admissible} if it is compatible with the operadic product in the following sense: if $m_1$,$m_2$ are monomials in $\P(V)([\zero, G])$ and $m'_1$, $m'_2$ monomials in $\P(V)([G, \un])$ for some $G \in P \in \P$ then $m_1 \trianglelefteq m_2, m'_1 \trianglelefteq m'_2$ implies $\mu_G(m_1, m'_1) \trianglelefteq \mu_G(m_2, m'_2)$. 
\end{madef}
\begin{ex}\label{exadmord}

Let $P$ be a poset in $\P$, let $\prec$ be some total order on $P$ (which may not have anything to do with the already existing order on $P$) and let $\propto$ be some total order on $\bigsqcup_{I \textrm{ interval of } P} B(I)$. One can define an admissible order $\trianglelefteq$ on monomials over intervals of $P$ as follow: 
\begin{multline*}
    \mu_{G_1, ..., G_n}(e_1, ... , e_n) \trianglelefteq \mu_{G'_1, ..., G'_{n'}}(e'_1, ... , e'_{n'}) \textrm{ if } G_1 \prec G'_1,  \\ \textrm{ or if } G_1 = G'_1 \textrm{ and } e_1 \propto e'_1, \\ \textrm{ or if } e_1 = e'_1, G_1 = G'_1, \textrm{ and } \mu_{G_2, ..., G_n}(e_2, ... , e_n) \trianglelefteq \mu_{G'_2, ..., G'_{n'}}(e'_2, ... , e'_{n'}). 
\end{multline*}
In other words we first compare the elements at the bottom of the chain, then the decorations at the bottom of the chain, and if they are equal we go up the monomial.
\end{ex}
Assume we have chosen an admissible ordering on monomials of $\P(V)$.
\begin{madef}[Leading term]

For any \red{nonzero} element $\alpha$ in $\P(V)(P)$ for some $P$ in $\nosP$, the \textit{leading term} of $\alpha$, denoted $\lt(\alpha)$, is the term $\lambda_{\alpha, m}m$ where $m$ is the greatest monomial with nonzero coefficient in $\alpha$ and $\lambda_{\alpha, m}$ is the coefficient of $m$ in $\alpha$. 
\end{madef}
We are finally ready for the main definition of this section. 
\begin{madef}[Gröbner basis]

Let $I$ be an ideal of $\P(V)$. A \textit{Gröbner basis} of $I$ is a subcollection $\G \subset I$ such that for any element $\alpha$ in $I$, the leading term of $\alpha$ is divisible by the leading term of some element in $\G$. 
\end{madef}
Note that Gröbner bases highly depend on the chosen basis of $V$, as well as the chosen admissible ordering on monomials. 
\begin{madef}[Normal monomial]

Let $\G$ be a subcollection of $\P(V)$. A \textit{normal monomial} of $\G$ is a monomial which is not divisible by the leading term of some element in $\G$. 
\end{madef}
\begin{myrmk}\label{rmqnormgen}
Using the rewriting procedure one can see that for any $\G$, if the admissible ordering on monomials is a well-order then the set of normal monomials of $\G$ linearly generates $\P(V)/\langle \G \rangle$.
\end{myrmk}
We have the following classical proposition. 
\begin{prop}\label{propmaingrob}
Let $\G$ be a subcollection of some ideal $I \subset \P(V)$. The subcollection $\G$ is a Gröbner basis of $I$ with respect to some order $\vartriangleleft$ if and only if the set of normal monomials of $\G$ with respect to $\vartriangleleft$ forms a basis of $\P(V)/I$.
\end{prop}
\subsection{Main results about Gröbner bases}
In this subsection we show that for particular choices of $\P$ the operads introduced in Section \ref{secdefex} admit quadratic Gröbner bases, which will be instrumental in subsequent sections. A key tool for finding Gröbner bases will be EL-labelings (see \cite{wachs_poset_2006}).
\begin{prop}\label{propgrobcom}
If every poset of $\P$ is EL-shellable then $\Com^{\P}$ admits the quadratic presentation $\Com^{\P}\simeq \P(V)/\langle R \rangle$ where $V$ is the $\P$-collection defined by 
\begin{equation*}
V(P) \coloneqq \left\{ 
\begin{array}{cl}
\Z & \textrm{ if } P \textrm{ has rank } 1, \\
\{0\} & \textrm{ otherwise,}
\end{array}
\right.
\end{equation*}
and $R$ is the collection of relations defined by 
\begin{equation*}
    R(P) \coloneqq \left \{
    \begin{array}{cl}
        \Z\langle \mu_{H, P}(1\otimes 1) - \mu_{H',P}(1\otimes 1), \, H, H' \in \iP\rangle & \textrm{ if } P \textrm{ has rank } 2,\\
    \{0\} & \mathrm{otherwise.} 
    \end{array}
    \right. 
\end{equation*}
Furthermore, $R$ is a Gröbner basis of $\langle R \rangle$. 
\end{prop}
\begin{proof}
In Example \ref{exPop} we have proved that the induced morphism $\P(V) \rightarrow \Com^{\P}$ is surjective and that $R$ is a subcollection of its kernel. Let $P$ be some poset in $\P$ and $\lambda: \E(P) \rightarrow P$ some edge-labelling of $P$. Monomials in $P$ can be identified with maximal chains in $P$. We order maximal chains in $P$ by lexicographic order on the corresponding words given by $\lambda$. By the definition of an edge-labelling, there is only one maximal chain whose associated word is increasing. All the other words contain a sequence of covering relations $G < G' < G''$ such that we have $\lambda(G<G') > \lambda(G'<G'')$. By the definition of an edge-labelling again, the monomial $G < G' < G''$ in $[G, G'']$ is not minimal for our ordering on monomials, and therefore is divisible by the leading term of some element in $R$. In other words, there is only one monomial in $P$ which is \blue{a normal monomial of $R$}. By Remark \ref{rmqnormgen} and Proposition \ref{propmaingrob} this proves that the morphism $\P(V)/\langle R \rangle \rightarrow \Com^{\P}$ is an isomorphism and that $R$ is a Gröbner basis of $\langle R \rangle$. 
\end{proof}
\begin{prop}\label{propgroblie}
If every poset of $\P$ is EL-shellable then the quadratic relations of $\Lie^{\P}$ form a Gröbner basis. \red{For every $P \in \nosP$, the normal monomials in $\Lie^{\P}(P)$ associated to that Gröbner basis have cardinality $|\mu(P)|$.}
\end{prop}
We postpone the proof of this result to Subsection \ref{subseckosdual}.
We now come to a central result of this article. 
\begin{prop}\label{propgrobgerst}
Let $\GL$ be the collection of geometric lattices. The operad $\Gerst^{\GL}$ is isomorphic to the operad $\OSop^{\vee}$ (see Example \ref{exPop}) and it admits a quadratic Gröbner basis. 
\end{prop}
We are specifically interested in the quadraticity of the Gröbner basis, as it will later imply that $\Gerst^{\GL}$ is Koszul (see Definition \ref{defkosbar}). 
\begin{proof}
This is where the combinatorics of geometric lattices shines through. We have a morphism of $\GL$-operads $\phi: \GL(V) \rightarrow \OSop^{\vee}$, with $V$ the generating collection of $\GGL$ (see Definition \ref{defgerst}), which is induced by the morphism 
\begin{equation*}
\begin{array}{ccc}
     V(\{\zero < \un\}) & \rightarrow &  (\OS(\{\zero < \un\}))^{\vee}  \\
       L &\rightarrow & e_{\un}^{\star} \\
       C &\rightarrow & 1^{\star}. 
\end{array}
\end{equation*}
Let us check that this morphism sends the ideal of relations $\langle R \rangle$ defining $\GGL$ to $0$. Let $\L$ be some geometric lattice of rank $2$ and $H, H'$ two atoms of $\L$. By definition of $\Phi$ we have 
\begin{equation*}
    \Phi(\mu_H(C\otimes C) - \mu_{H'}(C\otimes C)) = (1^{\star}\otimes 1^{\star})\circ \Delta_H - (1^{\star}\otimes 1^{\star})\circ \Delta_{H'}.
\end{equation*}
The two linear forms on the right hand side are equal (sending the unit of $\OS(\L)$ to $1$ and everything in strictly positive grading to $0$) and so the difference is $0$ as expected. On the other hand we have 
\begin{equation*}
    \Phi\left( \sum_{H \in \L^{\circ}} \mu_H(L \otimes L)\right) = \sum_{H \in \L^{\circ}}(e_{\un}^{\star}\otimes e_{\un}^{\star})\circ \Delta_H.
\end{equation*}
The linear form on the right is zero in grading other than $2$. Besides, it sends any element of the form $e_{H_1}e_{H_2}$ to $0$ (every term in the sum is zero except for $H= H_1$ and $H = H_2$ and those two terms cancel out) and therefore it is $0$ everywhere. Finally, we have 
\begin{equation*}
    \Phi\left( \mu_H(C\otimes L) - \sum_{H'\neq H}\mu_{H'}(L\otimes C)\right)  = (1^{\star}\otimes e_{\un}^{\star})\circ \Delta_H -  \sum_{H'\neq H} (e_{\un}^{\star} \otimes 1^{\star})\circ \Delta_{H'}.
\end{equation*}
and again the linear form on the right is $0$ everywhere (on $e_H$ every term is $0$ and on $e_{H'} \,(H' \neq H)$ we get $1-1 = 0$). Consequently, the morphism $\Phi$ induces a morphism on the quotient $\Phi: \GL(V)/\langle R \rangle \eqqcolon \GGL \rightarrow \OSop^{\vee}$.

Our next step is to prove that this morphism is surjective. This amounts to proving that for every $\L \in \GL$ the map 
\begin{equation*}
\begin{array}{ccc}
\OS(\L) &\rightarrow &\K\langle \textrm{ monomials of } \GL(V) \, \rangle \\
\alpha &\rightarrow & (\Phi(m)(\alpha))_m
\end{array}
\end{equation*}
is injective. Let $\vartriangleleft$ be a linear order on the atoms of $\L$. By \cite{Jambu_nbcbasis} any element $\alpha\in \OS(\L)$ can be uniquely written as a sum $\alpha = \sum_{B \textrm{ nbc-basis w.r.t. } \vartriangleleft} \lambda_B e_B$ (see \cite{Yuzvinsky_2001} for the definition of nbc-basis). We must construct monomials in $\GL(V)$ that will help us recover the coefficients $\lambda_B$ when applied to $\alpha$. Let $B = \{H_1 \vartriangleright \cdots \vartriangleright H_n \}$ be an nbc-basis. Let us denote by $c$ the maximal chain of $[\zero, H_1\vee \cdots \vee H_n]$
\begin{equation*}
    c\coloneqq \zero < H_1 < H_1\vee H_2 < \cdots < H_1\vee \cdots \vee H_n,
\end{equation*}
and $c'$ any maximal chain of $[H_1\vee \cdots\vee H_n, \un]$. Let us denote by $m$ the monomial 
\begin{equation*}
    m \coloneqq \mu_{\bigvee_i H_i}(\mu_c(L\otimes \cdots \otimes L)\otimes \mu_{c'}(C\otimes \cdots \otimes C)).
\end{equation*}
We will prove by induction on $n$ that we have the equality $\Phi(m)(\alpha) = \lambda_B$,
which implies the desired injectivity. One can check that we have $\Phi(m)(e_B) = 1$. Let $B'$ be an nbc-basis such that we have $\Phi(m)(e_{B'}) \neq 0$. We must prove the equality $B' = B$. By the definition of $m$ and the operadic product on $\OSop^{\vee}$ we readily have the inclusion $\{H_1\} \subset  B'$. Assume by induction that we have $
\{H_1, ... , H_k\} \subset  B'$ for some $k \leq n$. The inequality $ \Phi(m)(e_{B'}) \neq 0$ implies that there exists some element $H'_{k+1}$ in $B'$ such that we have $H_1 \vee ... \vee H_k \vee H_{k+1} = H_1\vee ... \vee H_k \vee H'_{k+1}$. If $H'_{k+1}$ is equal to $H_{k+1}$ then our induction step is complete. Otherwise there exists a circuit (see \cite{welsh_matroid_1976} for the definition) of $\L$ of the form $C = \{H_{i_1}, ..., H_{i_p}, H_{k+1}, H'_{k+1} \}$, for some indexes $i_1,..., i_p$ less than $k$. If $H_{k+1} \vartriangleright H'_{k+1}$ then the nbc-basis $B$ contains the broken circuit $C\setminus\{H'_{k+1}\}$ which is a contradiction. If not, the nbc-basis $B'$ contains the broken circuit $C\setminus \{H_{k+1}\}$ which is also a contradiction. This finishes the induction step, which proves that we necessarily have $B \subset B'$. By looking at the grading the cardinality of $B$ and $B'$ must be equal, and therefore we must have $B = B'$. This completes the proof of the surjectivity of~$\Phi$. \\

What is left is to prove the injectivity of $\Phi$ and finally find the claimed Gröbner basis of $\GGL$. As in the proof of Proposition \ref{propgrobcom} we will achieve both those goals at once. Let $\L$ be a geometric lattice and $\vartriangleleft$ a linear order on the set of atoms of $\L$. By \cite{Bjorner_1980} this order induces an EL-labeling $\lambda$ on~$\L$. Consider any admissible order on monomials $\vartriangleleft$ such that for any interval $[G, G']$ of rank $2$ in $\L$ we have the inequalities between monomials 
\begin{equation*}
    \mu_{H}(L\otimes C) \vartriangleleft \mu_{H'}(C\otimes L) \, \, \, \, \forall H, H' \in [G,G']^{\circ}, 
\end{equation*}
\begin{multline*}
    \mu_{H}(C\otimes C) \vartriangleleft \mu_{H'}(C\otimes C) \, \, \, \, \textrm{ for all } H, H' \in [G,G']^{\circ} \textrm{ s.t. } \lambda(G<H)\lambda(H<G') <_{lex} \lambda(G<H')\lambda(H'<G'),
\end{multline*}
and
\begin{multline*}
    \mu_{H}(L\otimes L) \vartriangleleft \mu_{H'}(L\otimes L) \, \, \, \, \textrm{ for all } H, H' \in [G,G']^{\circ} \textrm{ s.t. } \lambda(G<H)\lambda(H<G') >_{lex} \lambda(G<H')\lambda(H'<G') 
\end{multline*}
The corresponding normal monomials of $R$ are monomials of the form $\mu_G(m, m')$ where $G$ is some element in $\L$, $m$ is a normal monomial of $\Lie([\zero, G])$ for the order introduced in the proof of Proposition \ref{propgroblie} and $m'$ is a normal monomial of $\Com([G, \un])$ for the order introduced in the proof of Proposition \ref{propgrobcom}. By Propositions \ref{propgrobcom} and \ref{propgroblie} the cardinality of those monomials is $\sum_{G\in \L}|\mu([\zero, G])|$, which is also the dimension of $\OS(\L)^{\vee}$. By the surjectivity of $\Phi$, Remark \ref{rmqnormgen}, and Proposition \ref{propmaingrob} we get that $\Phi$ is an isomorphism and $R$ is a (quadratic) Gröbner basis of $\langle R \rangle$. \\
\end{proof}
We have the following corollary.
\begin{coro}\label{propdimgerst}
Let $\L$ be a geometric lattice. The Poincaré series of $\Gerst(\L)$, with grading given by the commutative generator, is $\sum_{F\in \L}|\mu(\zero, F)| t^{\rk [F, \un]}$.
\end{coro}

\section{Koszulness of $\P$-operads}\label{seckos}
For associative algebras and classical operads, Koszulness is a natural extension of quadraticity. In plain English, we say that a graded associative algebra/operad is Koszul if it is generated by elements of grading 1, relations between elements of grading 1 are generated by elements of grading 2, relations between relations of grading 2 are generated by elements of grading 3 and so on. For classical operads there are two main ways to formalize this, using either Koszul complexes or bar constructions (see \cite{LV_2012} for more details). In this section we will see that both  those approaches can be extended to $\P$-operads, with one extra assumption on $\P$: in the rest of this section we will assume that every poset of our collection $\P$ is \red{graded, meaning that for every $F \leq G$, all the maximal chains from $F$ to $G$ have the same length.} \red{In this section, all our operads will be defined in the monoidal category $\C$ of vector spaces over some field $\K$ with usual tensor product. In particular, by $\Lie$ and $\Com$ we will mean the tensorisation over $\K$ of the operads defined in Section \ref{secdefex}.}

\subsection{Bar construction}
\blue{For an associative algebra $A$, the bar construction of $A$ is the differential graded coalgebra defined as the cofree coalgebra over $A$, with differential given by the alternating sum of products of two elements (see \cite{polishchuk_quadratic_2005} Chapter 1). An algebra $A$ is called Koszul if the cohomology of its bar construction is concentrated in degree $0$.} In this subsection we define an analogous bar construction of $\P$-operads which leads to a notion of Koszulness in that context. Let us denote by $\dg-\C$ the monoidal category of complexes in $\C$. A $\P$-operad in $\dg-\C$ will be referred to as a $\dg$-$\P$-operad in $\C$. Note that by the Künneth formula the homology of a $\dg$-$\P$-operad is naturally a $\P$-operad. 
\begin{madef}[Bar construction]\label{defbarcons}
\phantom\qed Let $\O = (O, \mu)$  be a $\P$-operad. The bar construction of $\O$, denoted $\Bar(\O)$, is the $\dg$-$\P$-cooperad $(\P(O), \Delta, \d)$ where the operadic coproducts $\Delta_G \,  (G \in \iP \in \nosP)$  are defined by sending components $O([\zero,G_1])\otimes \cdots \otimes O([G_n, \un]) \subset \P(O)(P)$ to the same component viewed in $\P(O)([\zero, G])\otimes \P(O)([G, \un])$ (via the identity) if $G$ is one of the $G_i$'s, and sending the component to $0$ otherwise. The differentials $\d_P: \P(O)(P) \rightarrow \P(O)(P)\,(P\in \nosP)$ are defined on  components $O([\zero = G_0, G_1])\otimes \cdots \otimes O([G_n, G_{n+1} = \un])$ by 
\begin{equation*}
    \d (\alpha_0 \otimes \cdots \otimes \alpha_n) \coloneqq \sum_{i\leq n}(-1)^i \alpha_0 \otimes \cdots \otimes \alpha_{i-1} \otimes \mu_{G_{i+1}, [G_i, G_{i+2}]}(\alpha_i, \alpha_{i+1})\otimes \alpha_{i+2} \otimes \cdots \otimes \alpha_n,
\end{equation*}
which belongs to $\bigoplus_{i\leq n} O([\zero, G_1])\otimes \cdots \otimes O([G_i, G_{i+2}]) \otimes \cdots \otimes O([G_n, \un]) \subset \P(V)(P)$. 
\end{madef}
Those maps square to zero thanks to the associativity axiom \eqref{eqaxop} of $\P$-operads. If we assume furthermore that $\O$ is strictly positively graded then, thanks to the condition that the posets of $\P$ are graded one can put a cohomological degree on $\Bar(\O)$ by placing the summand $O_{i_0}([\zero, G_1])\otimes \cdots \otimes O_{i_n}([G_n,\un])$ in cohomological degree $i_0 + \cdots + i_n - (n+1)$. \blue{The differential $\d$ has degree $+1$ because the sum of the gradings $i_1 + \cdots + i_n$ does not change when applying $\d$ (because the operadic products are compatible with the grading) and the length $n$ of the chain decreases by one. We will denote the bar construction with this cohomological degree by $\Bar^{\bullet}(\O)$.} In the rest of this paper, every bar construction of a strictly positively graded operad will be given this degree.

\begin{madef}[Koszulness via bar construction]\label{defkosbar}
A strictly positively graded $\P$-operad $\O$ is said to be \textit{Koszul} if $\Bar^{\bullet}(\O)(P)$ has cohomology concentrated in degree $0$ for all $P$ in~$\nosP$.
\end{madef}
\begin{ex}\label{exordbar}
The operad $\Com^{\P}$ can be given a grading by placing $\Com^{\P}(P)$ in grading $\rk P$ for all $P \in \P$. One can see that for all $P$ in $\P$, the complex $\Bar^{\bullet}(\Com^{\P})(P)$ can be identified with the simplicial chain complex associated with the order complex of $P$ (put in cohomological degree convention) and therefore $\Com^{\P}$ is Koszul if and only if every poset of $\P$ is Cohen--Macaulay (see \cite{bjorner_introduction_1982} for the definition).
\end{ex}
\subsection{Koszul dual}\label{subseckosdual}
\blue{For an associative algebra $A$, the cohomology of the bar construction of $A$ is a coalgebra $A^{\kos}$ called the Koszul dual coalgebra of $A$, whose linear dual has a presentation given by the dual presentation of that of $A$ (see \cite{polishchuk_quadratic_2005} Section 2 Chapter 1). In this subsection we show that this is also the case in the context of $\P$-operads.}
\begin{madef}[Koszul dual (co)operad]

Let $\O$ be a \red{strictly positively graded} $\P$-operad. The \textit{Koszul dual cooperad} of $\O$, denoted $\O^{\text{!`}}$, is the $\P$-cooperad $H^{0}(\Bar^{\bullet}(\O))$. The \textit{Koszul dual operad} of $\O$, denoted  $\O^{!}$, is the $\P$- operad $(\O^{\kos})^{\vee}$. 
\end{madef}
Note that both make sense even if $\O$ is not Koszul. \red{The cooperad $\O^{\kos}$ is graded by the grading on $\Bar(\O)$ induced by that of $\O$.}
\begin{madef}[Twisting morphism]
For any \red{strictly positively graded} $\P$-operad $\O$ we have a morphism of $\P$-collections $\tau: \O^{\text{!`}} \rightarrow \O$ called the \textit{twisting morphism} of $\O$ which is the identity in grading $1$ and $0$ elsewhere. 
\end{madef}
\begin{prop}\label{propkosdual}
Let $\O = \P(V)/\langle R \rangle$ be a quadratic operad with $V$ some $\P$-collection. We have an isomorphism of operads
\begin{equation*}
    \O^{!} \cong \P(V^{\vee})/\langle R^{\bot}\rangle,
\end{equation*}
where $R^{\bot} \subset \P_{2}(V^{\vee})$ denotes the orthogonal of $R$ (over each poset $P \in \P$) for the pairing induced by the isomorphism $\P_{2}(V^{\vee}) \simeq (\P_{2}(V))^{\vee}$. 
\end{prop}
\begin{proof}
One has a  morphism $\P(V^{\vee})/\langle R^{\bot}\rangle \rightarrow \O^{!}$ coming from the universal property of $\P(V^{\vee})/\langle R^{\bot}\rangle$. This morphism is clearly surjective as it is induced by $\P(V^{\vee})\rightarrow \O^{!}$ which can be factored $\P(V^{\vee}) \xrightarrow{\simeq} \Bar^0(\O)^{\vee} \twoheadrightarrow \O^{!}$. The kernel of $\P(V^{\vee}) \rightarrow \O^{!}$ is exactly $\langle R^{\bot}\rangle$ by the general equality $(E\cap F)^{\bot} = E^{\bot}+F^{\bot}$.
\end{proof}
\begin{ex}\label{exdualliecom}

When $\P$ is such that $\Com^{\P}$ is quadratic, one immediately gets $$(\Com^{\P})^{!} \simeq \Lie^{\P}.$$
\end{ex}
\begin{ex}\label{exdualgerst}

For any $\P$, let $V$ and $R$ denote respectively the generating $\P$-collection and the relations of $\Gerst^{\P}$ as in Definition \ref{defgerst}. One easily computes 
\begin{equation*}
    R^{\bot}(P) \coloneqq  \left\{ 
\begin{array}{cl}
\begin{multlined}
    \Q \langle \mu_{H}(L^{\vee}\otimes L^{\vee}) - \mu_{H'}(L^{\vee}\otimes L^{\vee}), \, H, H' \in \iP,\\ \mu_H(L^{\vee}\otimes C^{\vee}) + \sum_{H' \neq H}\mu_{H'}(C^{\vee}\otimes L^{\vee}), \, H \in \iP,\\
    \sum_{H \in \iP} \mu_H(C^{\vee}\otimes C^{\vee}) \rangle
\end{multlined} & \textrm{ if } P \textrm{ has rank } 2, \\
\{0\} & \textrm{ otherwise.}
\end{array}
\right.
\end{equation*}
Note that $R^{\bot}$ is almost identical to $R$, with $L^{\vee}$ playing the role of $C$ and $C^{\vee}$ the role of $L$ (the only difference being the sign in the middle relation). When $\P$ is the collection of geometric lattices, by the same arguments as in the proof of Proposition \ref{propgrobgerst}  one can prove that $\Gerst^{!}$ is isomorphic to the linear dual of the cooperad $\textbf{tw}\OSop$ whose underlying $\GL$-collection is given by the Orlik--Solomon algebras and whose operadic coproducts are defined by  
\begin{equation*}
\Delta_G(e_H) \coloneqq \left \{
\begin{array}{cccc}
    e_H&\otimes& 1 & \textrm{ if } H \leq G,  \\
    - 1 &\otimes & e_{G\vee H} & \textrm{ otherwise.}
\end{array}
\right.
\end{equation*}
In particular $\Gerst^{!}$ has the same dimensions as $\Gerst$.

\end{ex}
We can now prove Proposition \ref{propgroblie} which we recall here.
\begin{prop}
If every poset of $\P$ is EL-shellable then the quadratic relations of $\Lie^{\P}$ form a Gröbner basis. \red{For every $P \in \nosP$, the normal monomials in $\Lie^{\P}(P)$ associated to that Gröbner basis have cardinality $|\mu(P)|$.}
\end{prop}
\begin{proof}
Consider $P$ a poset of $\P$ with $\EL$-labelling $\lambda$. Consider any order $\vartriangleright$ on monomials of $\Lie^{\P}(P)$ satisfying $\mu_{G, [F, H]}(1\otimes 1) \vartriangleright \mu_{G', [F, H]}(1\otimes 1)$ for any elements $G,G'$ in some rank $2$ interval $[F,H]$ of $P$ such that the word $\lambda(F,G)\lambda(G,H)$ is lexicographically before the word $\lambda(F,G')\lambda(G', H)$. The leading terms of quadratic relations for this order are quadratic monomials $\mu_{G, [F, H]}(1\otimes 1)$ where $F<G<H$ is the unique maximal chain with increasing labels in $[F,H]$. Normal monomials in $P$ with respect to this order are in bijection with maximal chains in $P$ with decreasing labels. By \cite{bjorner_wachs_1996} \red{Theorem 5.9} the cardinality of this set of normal monomials is the dimension of $H^{\rk P - 1}(P)$, which is the dimension of $\Lie^{\P}$ by Proposition \ref{propkosdual}, Example \ref{exordbar}, and Example \ref{exdualliecom}. By Proposition \ref{propmaingrob} this proves that the set of quadratic relations of $\Lie^{\P}$ forms a Gröbner basis with respect to $\vartriangleright$. \red{Besides, it is well known that for EL-shellable posets (or more generally for Cohen--Macaulay posets) the dimension of $H^{\rk P -1}(P)$ is $|\mu(P)|$ (see \cite{Hall_1936}) which gives the second statement.} 
\end{proof}
\subsection{Koszul complexes}
\blue{An alternative way of defining the Koszulness of an associative algebra $A$ is by asking the acyclicity of a complex $(A \otimes A^{\kos}, \d)$ called the Koszul complex of $A$ (see \cite{polishchuk_quadratic_2005} Chapter 2 Section 3). In this subsection we define an analogous notion for $\P$-operads and show that this gives an equivalent way of defining Koszulness in that context.}
\begin{madef}[Circle product on $\P$-collections]

Let $V$ and $W$ be two $\P$-collections. We define their product $V\circ W$ by 
\begin{equation*}
    V\circ W (P) \coloneqq \bigoplus_{G \in P}V([\zero, G])\otimes W([G, \un])
\end{equation*}
(with the convention that both $V(\{\star\})$ and $W(\{\star\})$ are equal to the monoidal unit of the chosen monoidal category). 
\end{madef}
\begin{myrmk}
This gives a monoidal structure on the category of $\P$-collections. One could have defined a $\P$-operad as a monoid in this monoidal category. 
\end{myrmk}
\begin{madef}[Koszul complex]\label{defkoscomplex}
Let $\O$ be a quadratic operad. The \textit{Koszul complex} $\Kos(\O)$ of $\O$ is the $\dg-\P$-collection $(O\circ O^{\text{!`}}, \d)$ with differentials \blue{$\d_P: O\circ O^{\kos}(P) \rightarrow O\circ O^{\kos}(P)$} defined on a component $O([\zero, G])\otimes O^{\kos}([G, \un]) \subset O\circ O^{\kos}(P)$  by
\begin{equation*}
    \d  \coloneqq \sum_{G' > G} (\mu_G\otimes \Id)\circ (\Id\otimes \tau \otimes \Id)\circ (\Id\otimes \Delta_{G'}),
\end{equation*}
with $\mu$ the operadic product of $\O$, $\Delta$ the operadic coproduct of $ \O^{\text{!`}}$ and $\tau$ the twisting morphism of $\O$. \blue{The Koszul complex of $\O$ has a degree given by the grading of $\O$, for which the differential $\d$ has degree $+1$.}
\end{madef} 
This definition is justified by the following lemma. 
\begin{mylemma}
    If $\O$ is quadratic then $\d$ squares to $0$. 
\end{mylemma}
\red{\begin{proof}
Denote $\O = \P(V)/\langle R \rangle $ with $R$ quadratic. We define $\tau^{(2)} := \sum_{G\in \iP}   \mu_G \circ(\tau \otimes \tau)\circ \Delta_G$, which is a morphism of $\P$-collections between $ \O^{\text{!`}}$ and $\O.$ By the associativity axiom one has 
\begin{align*}
    \d^2 = \sum_G  \mu_G\circ (\Id \otimes \tau^{(2)} \otimes \Id) \circ (\Id \otimes \Delta_G),
\end{align*}
so it is enough to prove that $\tau^{(2)}$ is zero. Since $\tau$ is zero in grading different from $1$, the morphism $\tau^{(2)}$ is zero in grading different from $2$. In grading equal $2$ by Proposition \ref{propkosdual} we can identify $\O^{\text{!`}}_2$ with $R$ and under this identification the morphism $\tau^{(2)}$ is the inclusion of $R$ in $\O$, which is zero. 
\end{proof}}
One has the following alternative definition of Koszulness. 
\begin{madef}[Koszulness via Koszul complexes]\label{defkoscx}

A quadratic operad $\O$ is said to be \textit{Koszul} if the Koszul complex $\Kos(\O) $ is acyclic. 
\end{madef}
\begin{prop}
    The two definitions of Koszulness Definition \ref{defkosbar} and Definition \ref{defkoscx} coincide.
\end{prop}
\begin{proof} 
Denote by $\overline{\Bar}(\O)$ the $\dg$-$\P$-collection 
$$\O^{\kos} \hookrightarrow \Bar^{0}(\O)\rightarrow \Bar^{1}(\O) \rightarrow \ldots.$$ Let us prove by induction on $n$ that $\overline{\Bar}(\O)$ is acyclic for every poset $P$ in $\P$ of rank less than $n$ if and only if the Koszul complex $\Kos(\O)(P)$ is acyclic for every poset $P$ in $\P$ of rank less than $n$. We introduce a decreasing filtration $\G^{\bullet}$ on $\Kos(\O)$ by setting 
\begin{equation*}
    \G^p(\Kos(\O)) = \bigoplus_{\substack{G \in \P \\ \rk G \geq p}}\O([\zero, G])\otimes \O^{\kos}([G, \un]).
\end{equation*}
The associated graded of this filtration is $\O\circ \O^{\kos}$ with zero differential, because there is no piece of the differential of $\Kos(\O)$ which does not strictly increase the rank of the interval given to $\O.$ Thus, as a graded vector space the first page $\Espec_1^p(\Kos(\O)) = H^{\bullet}(\gr^p \Kos(\O))$ is isomorphic to $\sum_{G \textrm{ of rank }p}\O([\zero, G])\otimes \O^{\text{!`}}([G, \un])$. Besides, the differential of $\Espec_1^p(\Kos(\O))$ is the piece of the differential of $\Kos(\O)$ which increases the rank of the interval given to $\O$ by exactly $1$. 

On the other hand, we can introduce a decreasing filtration $\F^{\bullet}$ on $\overline{\Bar}(\O)$ by setting 
\begin{equation*}
    \F^0\overline{\Bar}(\O) = \overline{\Bar}(\O), 
\end{equation*}
and 
\begin{equation*}
    \F^p\overline{\Bar}(\O) \coloneqq \bigoplus_{\substack{G_1 < \cdots < G_n \\ \rk G_1 \geq p}}\O([\zero, G_1])\otimes \cdots \otimes \O([G_n, \un])
\end{equation*}
for all $p \geq 1$. The associated graded of this filtration is 
$$\gr^0 \overline{\Bar}(\O) = \cdots \rightarrow 0\rightarrow \O^{\text{!`}} \rightarrow 0 \rightarrow \cdots, $$ and
$$\gr^p \overline{\Bar}(\O) \cong \sum_{G \textrm{ of rank } p} \O([\zero, G])\otimes \Bar^{\bullet}(\O)([G, \un])$$ for all $p \geq 1.$ If we assume that $\overline{\Bar}(\O)$ is acyclic for every poset of rank less than or equal to $n-1$, then for every poset of rank less than or equal to $n$ the first page $\Espec_1^p(\overline{\Bar}(\O) = H^{\bullet}(\gr^p \overline{\Bar}(\O))$ is isomorphic, as a graded vector space, to $\sum_{G \textrm{ of rank }p}\O([\zero, G])\otimes \O^{\text{!`}}([G, \un])$ (with the convention $\O([\zero, \zero]) = \K$). Under this identification, the differential $\d^1$ on a summand $\O([\zero, G])\otimes \O^{\kos}([G, \un])$ is $$\sum_{\substack{G'> G \\ \rk G' - \rk G = 1}} (\mu_G \otimes \Id)\circ (\Id \otimes \tau \otimes \Id)\circ (\Id \otimes \Delta_{G'}),$$ and so for all $p$ we have the isomorphism of dg complexes 
\begin{equation*}
    (\Espec_1^p(\overline{\Bar}(\O)), \d^1) \cong (\Espec_1^p(\Kos(\O)), \d^1).
\end{equation*}
More generally, one can check that we have this isomorphism for every page (the differential at the page $k$ being the piece of $\d$ which increases the rank of the interval given to $\O$ by exactly $k$). By a standard spectral sequence argument this proves that $\overline{\Bar}(\O)$ is acyclic for every poset of rank $\leq n$, if and only if $\Kos(\O)$ is acyclic for every poset of rank $\leq n$, which concludes the proof. 
\end{proof}

\subsection{Gröbner bases and Koszul duality}
For associative algebras and classical operads, we have the key proposition that admitting a quadratic Gröbner basis implies being Koszul (see \cite{Hoffbeck_2010} for operads for instance). In this subsection we prove that this is also true for $\P$-operads. This result will be our main tool for proving Koszulness of $\P$-operads. 

\begin{prop}\label{propgrobkos}
    Let $\O$ be a strictly positively graded $\P$-operad. If $\O$ is generated by elements of grading $1$ and admits a quadratic Gröbner basis then $\O$ is Koszul.
\end{prop}

\red{\begin{proof}
This is just an adaptation of the proof given in \cite{Hoffbeck_2010} to the setting of $\P$-operads. Let us denote $\O \cong \P(O_1)/\langle \G \rangle$ where $\G$ is a quadratic Gröbner basis of the ideal $\ker \P(O_1) \rightarrow \O$, for some choice $\B = \{e_i\}$ of a basis of $O_1$ and some choice $\vartriangleleft$ of an admissible well-order on monomials (see Section \ref{secgrobner} for the \blue{definitions}). We introduce a filtration $\F$ on $\Bar^{\bullet}(\O)$ indexed by monomials of $\O$, by setting 
\begin{equation*}
    \F^{m}(\Bar^{\bullet}(\O))  \coloneqq \K\langle \{ m_1\otimes \cdots \otimes m_n \, | \, m_i \textrm{ monomials such that } \mu(m_1 \otimes \cdots \otimes m_n) \trianglelefteq m \} \rangle
\end{equation*}
This is an increasing filtration compatible with the differential of the bar construction. Let us compute the associated graded $\gr^m \Bar^{\bullet}(\O)$ for some monomial $m = e_0\otimes \cdots \otimes e_n \in O_1([\zero, G_1])\otimes \cdots \otimes O_1([G_n, \un])$, with $G_1<\cdots < G_n$ some chain in some poset $P\in \P$. For $I$ a subset of consecutive elements of $\{0, \ldots, n\}$ we denote by $e_I$ the operadic product of the generators $e_i, i \in I$. For $\pi$ an ordered partition of $\{0, \ldots,n\}$ with blocks of consecutive elements $\pi_1,\ldots,\pi_k$ we denote by $m_{\pi}$ the element $e_{\pi_1}\otimes \cdots \otimes e_{\pi_k}$ of the bar construction. Let us denote by $\Adm(m)$ the set of indexes $i \leq n-1 $ such that $\mu(e_i\otimes e_{i+1})$ is a normal monomial in $\O$. By definition of $\F$ the associated graded $\gr^m\Bar^{\bullet}(\O)(P)$ is linearly generated by the elements $m_{\pi}$ where $\pi$ runs over the ordered partitions of $\{0, \ldots,n\}$ such that for every block $\pi_i\in \pi$ with more than two elements, the set $\pi_i\setminus \{\max \pi_i\}$ is included in $\Adm(m)$. For such partitions the elements $e_{\pi_i}$ have no quadratic non normal submonomials and so they are normal, by the quadraticity of $\G$. This means that the elements $m_{\pi}$ form a linear basis of $\gr^m\Bar^{\bullet}(\O)(P)$. As a consequence, the map $\pi \rightarrow \bigcup_i \pi_i\setminus \{\max \pi_i\}$ induces an isomorphism of $\dg$-complexes between $\gr^m\Bar^{\bullet}(\O)(P)$ and the augmentation of $C^{\bullet}(\Delta_{\Adm(m)})$ \red{(where $\Delta$ denotes the standard simplex)}. This latter complex has cohomology zero unless $\Adm(m)$ is empty. In this case $\gr^m(\Bar^{\bullet}(\O))$ is equal to $\K$ concentrated in degree $0$ (generated by $e_0\otimes \cdots \otimes e_n$). By a standard spectral sequence argument this implies that the cohomology of $\Bar^{\bullet}(\O)$ is concentrated in degree $0$, which concludes the proof.
\end{proof}}
By virtue of Proposition \ref{propgrobgerst} this implies the following. 
\begin{coro}\label{corokosgerst}
    The $\GL$-operad $\Gerst^{\GL}$ is Koszul.
\end{coro}

\section{Application to Kazhdan--Lusztig--Stanley theory}\label{secKLS}
In this section we explain how the constructions of the previous sections can be used to get a categorification of the Kazhdan--Lusztig polynomials of geometric lattices introduced in  \cite{elias_kazhdanlusztig_2016}.

\subsection{Reminders}
In this subsection we briefly outline the theory of Kazhdan--Lusztig--Stanley polynomials as introduced in \cite{Stanley_1992}. We refer to \cite{proudfoot_algebraic_2018} for more details. Let $P$ be a locally finite poset which is graded. Let $I_{\rk}(P)$ (resp. $I_{\rk/2}(P)$) be the subring of the incidence algebra $I(P) \coloneqq \prod_{G_1\leq G_2 \in P}\Z[X]$ which consists of elements $f$ such that $f_{G_1G_2}$ has degree \blue{less than or equal to} $\rk [G_1, G_2]$ for all $G_1 \leq G_2 \in P$ (resp. strictly less than $\rk [G_1, G_2]/2$). The subring $I_{\rk}(P)$ admits an involution $f\rightarrow \overline{f}$ defined by 
\begin{equation*}
    \overline{f}_{G_1G_2} = t^{\rk([G_1, G_2])}f_{G_1G_2}(t^{-1}).
\end{equation*}
We denote by $\delta$ the unit of $I(P)$, which is equal to $1$ on every interval which is a singleton, and $0$ elsewhere. 
\begin{madef}[P-kernel]

A $P$-\textit{kernel} $\kappa$ is an element of $I_{\rk}(P)$ \red{which is equal to $1$ on singleton intervals,} and which satisfies the equation $$\overline{\kappa} \kappa = \delta .$$ The following theorem is a consequence of \cite{Stanley_1992} Corollary 6.7.
\end{madef}
\begin{theo}\label{theoStanley}
Let $\kappa$ be \red{a $P$-kernel}. There exists a unique pair of element\red{s} $f, g \in I_{\rk/2}$ such that we have $\overline{f} = \kappa f$, $\overline{g} = g \kappa$, and $f_{GG} = g_{GG} = 1 \, \, \forall \, G \in P$.
\end{theo}
Stanley only considered the equation $\overline{g} = g \kappa$ which he called $\kappa$-acceptability. Following \cite{brenti_twisted_1999} we will call the polynomials $f$ and $g$ respectively the right and left KLS polynomials \blue{of $P$ with respect to $\kappa$}. If $\kappa$ is a $P$-kernel then $\overline{\kappa}$ is also a $P$-kernel whose right and left KLS polynomials are called inverse right and left KLS polynomials. Note that 
\begin{ex}
\begin{itemize}
    \item The characteristic polynomial of the intervals of $P$ is a $P$-kernel (\cite{Stanley_1992}, Example 6.8). If $P$ is a geometric lattice then the corresponding right KLS polynomial is the Kazhdan--Lusztig polynomial of $P$ introduced in \cite{elias_kazhdanlusztig_2016}.
    \item If $P$ is an Eulerian poset then the element $\kappa \in I(P)$ defined by $\kappa_{G_1G_2}(t) \coloneqq (t - 1)^{\rk [G_1, G_2]}$ for all $G_1 \leq G_2$ is a $P$-kernel (\cite{Stanley_1992}, Proposition 7.1). In the case where $P$ is the face lattice of a polytope $\Delta$, the corresponding left KLS polynomial is the $g$-polynomial of $\Delta$. 
    \item If $W$ is a Coxeter group with its Bruhat order, its $R$-polynomial (see \cite{bjorner_combinatorics_2005} Chapter 5) is a $W$-kernel. The corresponding right KLS polynomial is the classical Kazhdan--Lusztig polynomial defined in \cite{KL_1979}.
\end{itemize}
\end{ex}
\subsection{Categorification of KLS theory}
Recall from Corollary  \ref{propdimgerst} that for any geometric lattice $\L$ the Poincaré series of $\Gerst^{\GL}(\L)$ (with grading given by the commutative generator) is $\sum_{F \in \L}|\mu(\zero, F)|t^{\rk \L - \rk F}.$ Notice that this is essentially the characteristic polynomial $\chi_{\L} \coloneqq \sum_{F \in \L}\mu(\zero, F)t^{\rk \L - \rk F}$ of $\L$, up to an alternating sign. Recall from Corollary \ref{corokosgerst} that $\Gerst^{\GL}$ is Koszul. Using the definition of Koszulness by Koszul complexes (Definition \ref{defkoscx}) this means that the Koszul complex $\Kos(\Gerst)$ is acyclic, which implies that its graded Euler characteristic is zero. By the study of $\GGL^{!}$ carried out in Example \ref{exdualliecom} this exactly means that the convolution product $(\overline{\chi}\chi)_{\L}$ is zero for non trivial $\L$ (the bar coming from the fact that in $\Gerst^{!}$ the role of the commutative generator is played by $L^{\vee}$ instead of $C^{\vee}$). This recovers the fact that the characteristic polynomial is a $P$-kernel. Having a new proof of this somewhat elementary result is not interesting in itself, but it is our first hint of a connection between the theory of $\P$-operads developed in this article, and the theory of \blue{Kazhdan--Lusztig polynomials of geometric lattices}. In this section, we refine this connection to provide a categorification of the KLS polynomials themselves.\\

Let us go back to the defining equation for, say, right KLS polynomials:
\begin{equation*}
    \overline{f} = \kappa f.
\end{equation*}
The convolution product on the right contains the term $f_{\zero, \un}$. Putting this term on the left we get 
\begin{equation*}
    \overline{f}_{\zero, \un} - f_{\zero, \un} = \sum_{G > \zero}\kappa_{\zero, G}f_{G,\un}.
\end{equation*}
At this point one has to remember the very crucial fact that KLS polynomials have degree strictly less than half the rank of $P$. This means that $\overline{f}_{\zero, \un}$ and $f_{\zero, \un}$ are supported on different degrees and we can define $f$ as the truncation of $-\sum_{G > \zero}\kappa_{\zero, G}f_{G,\un}$ in degree strictly less than half the rank of $P$. \red{If we denote by $\tr$ this truncation operation and we iterate this formula we get 
\begin{equation*}
    f_{\zero, \un} = \sum_{\substack{n \in \N, \\\zero = G_0 < \cdots < G_n = \un}}(-1)^n \tr\big(\kappa_{\zero, G_1}\tr\big(\kappa_{G_1, G_2}\tr\big( \ldots \big)\ldots \big).
\end{equation*}}
If we imagine that $\kappa$ has been categorified by some operad $\O$ (maybe up to an alternating sign, as in the case of the characteristic polynomial described above) this strongly suggests that we should look for a categorification of $f$ in the bar construction $\Bar(\O)$ of $\O$ (see Definition \ref{defbarcons}). In fact, it directly suggests the following definitions.
\begin{madef}[KLS complexes]\label{defKLScomplex}

Let $\O_{(\bullet,\bullet)}$ be a \blue{positively} bigraded $\P$-operad for some fbs collection $\P$ of graded posets. We define the sub-complex $\RKLS_{\O}$ of $\Bar(\O)$ by declaring that a graded summand 
\begin{equation}\label{eqsumRKLS}
    \bigotimes_{0 \leq k\leq n-1}\O_{(i_k, j_k)}([G_k, G_{k+1}]) \subset \Bar(\O)(P)
\end{equation}
belongs to $\RKLS_{\O}(P)$ if we have $\sum_{p\geq q}i_p < \rk [G_q, \un]/2$ for all $q >0$. Similarly, we define the sub-complex $\LKLS_{\O}(P)$ by asking the condition $\sum_{p\leq q}i_p < \rk [\zero, G_{q+1}]/2$ for all $q < n-1$. Finally, we also define $\iRKLS_{\O}$ and $\iLKLS_{\O}$ by swapping the two gradings in the conditions defining $\RKLS_{\O}$ and $\LKLS_{\O}$ respectively. 
\end{madef}
If the operad $\O$ is clear from the context we will omit it.
\blue{\begin{rmq}
The differential of some element in some summand \eqref{eqsumRKLS} consists of a sum of elements each living in a summand of the form 
\begin{multline*}
     \O_{(i_0,j_0)}([G_0, G_1]) \otimes \cdots \otimes \O_{(i_{k-1}, j_{k-1})}([G_{k-1}, G_{k}])\otimes \O_{(i_{k}+i_{k+1}, j_k + j_{k+1})}([G_{k}, G_{k+2}])\otimes \\ \O_{(i_{k+2}, j_{k+2})}([G_{k+2}, G_{k+3}]) \otimes \cdots \otimes \O_{(i_{n-1}, j_{n-1})}([G_{n-1}, G_n])
\end{multline*}
for some $k\leq n-2$. If the summand \eqref{eqsumRKLS} belongs to $\RKLS$ (resp. $\LKLS$, $\iRKLS$, $\iLKLS$) then all those summands belong to $\RKLS$ (resp. $\LKLS$, $\iRKLS$, $\iLKLS$) and so $\RKLS$ (resp. $\LKLS$, $\iRKLS$, $\iLKLS$) is a subcomplex of $\Bar(\O)$.
\end{rmq}}
Those complexes have a bigrading induced by the bigrading of $\O$. \red{The part of the complex $\RKLS(P)$ (resp. $\LKLS(P)$) of first grading strictly less than $\rk \, P/2$ is meant to categorify the right (resp. left) KLS polynomial, and $\widehat{\RKLS}$, $\widehat{\LKLS}$ their inverse versions.} \\

Let us now turn our attention toward the case $\O = \GGL$, bigraded by the two generators $C$ and $L$. In this case $\GGL_{(p,q)}(\L)$ is non-trivial only if $p+q = \rk \L$ so in our notation we can forget one of the two gradings, say the one given by the Lie generator. 
\begin{nota*}
    We will refer to the grading given by the commutative generator as the weight.
\end{nota*} 
Let us describe the first few KLS complexes. We use the notation 
\begin{equation*}
    \begin{array}{ccc}
     C^{i_{n-1}}L^{j_{n-1}} & & \\ 
     \cdot \cdot \cdot &  \coloneqq &  \bigoplus_{\zero = G_0< \cdots < G_n = \un }\bigotimes_{0 \leq k\leq n-1}\GGL_{(i_k, j_k)}([G_k, G_{k+1}]). \\
     C^{i_{0}}L^{j_{0}} & &
     \end{array}
\end{equation*}
Such a summand of $\Bar(\GGL)(P)$ lives in cohomological degree the rank of $P$ minus the number of storeys. It belongs to $\RKLS(P)$ if for each proper upper interval, strictly less than half of the letters are $C$'s.
\begin{ex}\label{exKLScx}

If $\L$ is of rank $1$ we have 
\begin{equation*}
    \RKLS_{(0)}(\L) = \iRKLS_{(1)}(\L) = L, \RKLS_{(1)}(\L) = \iRKLS_{(0)}(\L) = C.
\end{equation*}
and the left KLS complexes are the same as their right counterpart. 
If $\L$ is of rank $2$ we have 
\begin{equation*}
\begin{tikzcd}[column sep = tiny, row sep = tiny]
    \RKLS_{(0)}(\L) = \begin{array}{c} L\\ L \end{array} \longrightarrow L^2,
\end{tikzcd}
\begin{tikzcd}[column sep = tiny]
    \RKLS_{(1)}(\L) = \begin{array}{c} L\\ C \end{array} \longrightarrow  CL,
\end{tikzcd}
\RKLS_{(2)}(\L) = C^2.
\end{equation*}
From $\RKLS$ one can get the left KLS complexes by flipping top and bottom, and the inverse KLS complexes by exchanging $L$ and $C$. If $\L$ is of rank $3$ we have 
\begin{equation*}
\RKLS_{(0)}(\L) = 
\begin{tikzcd}[row sep = tiny]
& \begin{array}{l} L^2\\ L \end{array} \arrow[rd] & \\
\begin{array}{l} L\\ L \\ L \end{array} \arrow[ru] \arrow[rd] & & L^3\\ 
     & \begin{array}{l} L\\ L^2 \end{array}  \ar[ru] &
\end{tikzcd}
\end{equation*}
which is isomorphic to $\Bar^{\bullet}(\Lie)(\L)$, 
\begin{equation*}
\RKLS_{(1)}(\L) = 
\begin{tikzcd}[row sep = tiny]
& \begin{array}{l} L^2\\ C \end{array} \arrow[rd] & \\
\begin{array}{l} L\\ L \\ C \end{array} \arrow[ru] \arrow[rd] & & CL^2,\\ 
     & \begin{array}{l} L\\ CL \end{array}  \ar[ru] & 
\end{tikzcd}
\end{equation*}
\begin{equation*}
\RKLS_{(2)}(\L) =
\begin{tikzcd}
   \begin{array}{l} L\\ C^2 \end{array} \arrow[r] & C^2L,
\end{tikzcd}
\end{equation*}
and finally
\begin{equation*}
    \RKLS_{(3)}(\L) = C^3. 
\end{equation*}
The other KLS complexes can be obtained from $\RKLS$ using the same transformations as in rank $2$ (exchanging $L$ and $C$ and/or top and bottom). Let us finish with a KLS complex in rank $4$. If $\L$ is of rank $4$ we have
\begin{equation*}
    \RKLS_{(2)}(\L) = 
\begin{tikzcd}[row sep = tiny]
& \begin{array}{l} L^2\\ C \\ C  \end{array} \ar[r] \ar[rd] & \begin{array}{l} CL^2\\ C \end{array} \ar[rd] & \\
\begin{array}{l} L\\ L \\ C \\ C   \end{array} \ar[ur] \ar[r] \ar[dr] &  \begin{array}{l} L \\ CL \\ C \end{array} \ar[ur] \ar[dr] &  \begin{array}{l} L^2\\ C^2 \end{array} \ar[r] & C^2L^2. \\
&  \begin{array}{l} L \\ L \\  C^2\end{array} \ar[r] \ar[ur] &  \begin{array}{l} L\\ C^2L \end{array} \ar[ur]
\end{tikzcd}
\end{equation*}
\end{ex}
By construction of the KLS complexes we have the following lemma.
\begin{mylemma}\label{lemmaeulerchar}
Let $\L$ be a geometric lattice, $P_{\L}(t)$ its Kazhdan--Lusztig polynomial as defined in \cite{elias_kazhdanlusztig_2016} (that is, the left KLS polynomial with respect to the $\L$-kernel $\chi_{\L}$) and $Q_{\L}(t)$ its inverse Kazhdan--Lusztig polynomial (right KLS polynomial with respect to the $\L$-kernel $\overline{\chi_{\L}}$). We have the identities
\begin{equation*}
    P_{\L}(t) = \sum_{i< \rk \L /2}(-1)^{i}\chi_{\Eulerchar}(\RKLS_{(i)})t^i, \, \,  Q_{\L}(t) = \sum_{i< \rk \L /2}(-1)^{i}\chi_{\Eulerchar}(\iLKLS_{(i)})t^i,
\end{equation*}
where $\chi_{\Eulerchar}$ denotes the Euler characteristic. 
\end{mylemma}
\blue{\begin{proof}
Let us denote $$\chi^{+}_{\L}(t) \coloneqq \sum_{G \in \L}|\mu([\zero, G])|t^{\rk [G, \un]} = \Hilb(\Gerst(\L))(t) = (-1)^{\rk \, \L}\chi_{\L}(-t).$$ 
Note that for any chain $\zero = G_0 < \cdots < G_n = \un$ in $\L$, the Hilbert series of the subspace of $\RKLS(\L)$ consisting of summands with underlying chain $\zero = G_0 < \cdots < G_n = \un$ is $\tr\big(\chi^+_{\zero, G_1}(t)\tr\big(\chi^+_{G_1, G_2}(t) \tr\big(\ldots\big)\ldots\big)$. Then, for every geometric lattice $\L$ we have 
\begin{align*}
    P_{\L}(t) &= \sum_{\zero = G_0 < \cdots < G_n = \un}(-1)^n \tr\big(\chi_{\zero, G_1}(t)\tr\big(\chi_{G_1, G_2}(t) \tr\big(\ldots\big)\ldots\big) \\
    &= \sum_{\zero = G_0 < \cdots <G_n = \un}(-1)^{\rk \L - n }\tr\big(\chi^+_{\zero, G_1}(-t)\tr\big(\chi^+_{G_1, G_2}(-t) \tr\big(\ldots\big)\ldots\big)  \\
    &= \sum_{k}(-1)^k\Hilb(\RKLS^k(\L))(-t) \\
    &= \chi_{\Eulerchar}(\RKLS)(\L)(-t) \\
    &= \sum_{i}(-1)^i\chi_{\Eulerchar}(\RKLS_{(i)})t^i.
\end{align*}
The case of left, or inverse KLS polynomials are similar. 
\end{proof}}
We can now state the main result of this article. 
\begin{theo}\label{theomainhomo}
Let $\P$ be a collection of finite bounded graded posets stable under taking closed intervals, such that the operad $\Gerst^{\P}$ is Koszul. Let $P$ be a poset in $\P$ of rank $k$. 
\begin{enumerate}[i)]
    \item For $i<k/2$ the cohomology of $\RKLS_{(i)}(P)$, $\LKLS_{(i)}(P)$, $\iRKLS_{(i)}(P)$ and $\iLKLS_{(i)}(P)$ is concentrated in degree $i$. For $i> k/2$ the cohomology of $\RKLS_{(i)}(P)$, $\LKLS_{(i)}(P)$, $\iRKLS_{(i)}(P)$ and $\iLKLS_{(i)}(P)$ is concentrated in degree $i-1$. 
    \item If $k$ is even, the complexes $\RKLS_{(\frac{k}{2})}(P)$, $\LKLS_{(\frac{k}{2})}(P)$, $\iRKLS_{(\frac{k}{2})}(P)$ and $\iLKLS_{(\frac{k}{2})}(P)$ are acylic. 
\end{enumerate}
\end{theo}
Of course we only have one example of a collection $\P$ such that $\Gerst^{\P}$ is Koszul, namely the collection $\GL$ of geometric lattices (Corollary \ref{corokosgerst}), and it wouldn't be surprising to the author if every such collection $\P$ was contained in $\GL$. We state the theorem with this degree of generality to emphasize the fact that the combinatorics of geometric lattices needed in this article is completely contained in Corollary \ref{corokosgerst} (and by extension in Proposition \ref{propgrobgerst}, which roughly boils down to the shellability of geometric lattices and the existence of a cooperadic structure on Orlik--Solomon algebras). 

Statement i) of Theorem \ref{theomainhomo}, Corollary \ref{corokosgerst} and Lemma \ref{lemmaeulerchar} give a new proof of the following result. 
\begin{coro}[\cite{braden2023singular}]
Let $\L$ be a geometric lattice. The polynomials $P_{\L}$ and $Q_{\L}$ have positive coefficients. 
\end{coro}
\begin{proof}[Proof of Theorem \ref{theomainhomo}]
We will prove the two statements together by induction on $\rk P$. The base cases are $i=0$, $i= \rk P$ in any rank for Statement i), and $i=1$ in rank $2$ for Statement ii). For the base cases of Statement i) the KLS complexes are isomorphic to either $\Bar(\Lie)$, $\Bar(\Com)$, $C^{\rk P}$ (one term complex) or $L^{\rk P}$ (one term complex). By the Koszulness of $\Lie$ and $\Com$ those complexes have cohomology concentrated in the expected degree. For the base case of Statement ii) we must prove that the morphisms 
\begin{equation*}
\begin{tikzcd}
   \begin{array}{c} L\\ C \end{array} \arrow[r] & CL, 
\end{tikzcd}
\textrm{ and }
\begin{tikzcd}
   \begin{array}{c} C\\ L \end{array} \arrow[r] & CL, 
\end{tikzcd}
\end{equation*}
are isomorphisms for any poset of $\P$ of rank $2$. The graded summand $CL$ can be naturally identified with $\Q\langle \textrm{Atoms of } P\rangle$ by sending an atom $H$ to the monomial $\mu_{H}(L, C)$, and so can $\begin{array}{c} C\\ L \end{array}$ and $\begin{array}{c} L\\ C \end{array}$, by definition. Under those identifications the two morphisms above are respectively 
\begin{equation*}
    \begin{array}{ccc}
         \Q \langle \textrm{Atoms of } P \rangle & \rightarrow &  \Q \langle \textrm{Atoms of } P \rangle\\
         H & \rightarrow & \sum_{H' \neq H} H',  
    \end{array}
\end{equation*}
and
\begin{equation*}
    \begin{array}{ccc}
         \Q \langle \textrm{Atoms of } P \rangle & \rightarrow &  \Q \langle \textrm{Atoms of } P \rangle\\
         H & \rightarrow & H.  
    \end{array}
\end{equation*}
Since we are working over $\Q$, both those morphisms are isomorphisms. 

The general strategy for the induction step of Statement i) is to find an increasing filtration $\F^{\bullet}$ of $\Bar(\GGL)$ such that we have $$\RKLS_{(i)} = \F_{(i)}^{p}\Bar(\GGL)$$ for some $p$, such that every associated graded complex in filtration grading higher than $p$ has cohomology concentrated in degree less than $i-1$ if $i < \rk /2$ (resp. $i-2$ if $i > \rk /2$) and such that every associated graded complex in filtration grading lower than $p$ has cohomology concentrated in degree greater than $i$ if $i < \rk /2$ (resp. $i-1$ if $i> \rk /2)$. If such a filtration exists then by a standard spectral sequence argument the cohomology of $\RKLS_{(i)}$ is concentrated in degrees higher than $i$ (resp. $i-1$) and the cohomology of the quotient complex $\Bar(\GGL)_{(i)}/\RKLS_{(i)}$ is concentrated in degrees lower than $i-1$ (resp. $i-2$) which implies by the associated long exact sequence and the acyclicity of $\Bar(\GGL)$ in strictly positive degrees that the cohomology of $\RKLS_{(i)}$ is concentrated in degree $i$ (resp. $i-1$) as claimed. The filtration we will define has a natural interpretation in terms of lattice paths. 
\begin{madef}[Lattice path]

A \textit{lattice path} is a map from a finite subset of $\Z$ to $\Z$.
\end{madef}

Each graded summand 
\begin{equation*}
     \begin{array}{ccc}
     C^{i_{n-1}}L^{j_{n-1}} & & \\ 
     \cdot \cdot \cdot &  \coloneqq &  \bigoplus_{\zero = G_0 < \cdots < G_n = \un}\bigotimes_{0 \leq k\leq n-1}\GGL_{(i_{k}, j_{k})}([G_{k}, G_{k+1}]) \\
     C^{i_{0}}L^{j_{0}} & &

    \end{array}
\end{equation*}
of $\Bar(\GGL)(P)$ has an associated lattice path $\varphi$ defined inductively by 
\begin{equation*}
\left\{ \begin{array}{ccc}
  \varphi(0) & = & 0 \\
  \varphi(\rk P - \rk G_{q}) & = & \varphi(\rk P - \rk G_{q+1}) + j_{q} - i_{q}.
\end{array}
\right.
\end{equation*}
In other words, we start at zero and we read the symbols of 
\begin{equation}\label{eqgradedsummand}
\begin{array}{c}
    C^{i_{n-1}}L^{j_{n-1}} \\
    \cdots  \\
     C^{i_{0}}L^{j_{0}} \\
\end{array}
\end{equation}
from top to bottom. When we see an $L$ we move in direction $(1,1)$ and when we see a $C$ we move in direction $(1,-1)$. Finally, we forget the values in between the storeys of \eqref{eqgradedsummand}. Here are some examples below. 

\begin{ex}

\begin{equation*}
\begin{array}{cc} L  \\ L  \\ C \end{array}  \rightarrow
\begin{tikzpicture}[scale=0.60, baseline=14px]
    \draw (0,0) grid (3,2);
    \draw[ultra thick,red]
 (0,0) -- (1,1) -- (2,2) -- (3,1);
    \fill[red] (0,0) circle (5pt)
               (1,1) circle (5pt)
               (2,2) circle (5pt)
               (3,1) circle (5pt);
\end{tikzpicture}
\end{equation*}
\begin{equation*}
\begin{array}{cc} L  \\ C  \\ L \end{array}  \rightarrow
\begin{tikzpicture}[scale=0.60, baseline=14px]
    \draw (0,0) grid (3,2);
    \draw[ultra thick,red]
 (0,0) -- (1,1) -- (2,0) -- (3,1);
    \fill[red] (0,0) circle (5pt)
               (1,1) circle (5pt)
               (2,0) circle (5pt)
               (3,1) circle (5pt);
\end{tikzpicture}
\end{equation*}
\begin{equation*}
\begin{array}{cc} L^2  \\ C \end{array}  \rightarrow
\begin{tikzpicture}[scale=0.60, baseline=14px]
    \draw (0,0) grid (3,2);
    \draw[ultra thick,red]
 (0,0) -- (2,2) -- (3,1);
    \fill[red] (0,0) circle (5pt)
               (2,2) circle (5pt)
               (3,1) circle (5pt);
\end{tikzpicture}
\end{equation*}
\begin{equation*}
\begin{array}{cc} CL  \\ L \end{array}  \rightarrow
\begin{tikzpicture}[scale=0.60, baseline=14px]
    \draw (0,0) grid (3,2);
    \draw[ultra thick,red]
 (0,0) -- (2,0) -- (3,1);
    \fill[red] (0,0) circle (5pt)
               (2,0) circle (5pt)
               (3,1) circle (5pt);
\end{tikzpicture}
\end{equation*}
\end{ex}

\begin{myrmk}\label{rmkweightpath}
\phantom\qed
This assignment is obviously injective and in particular one can retrieve the numerical parameters of a graded summand from the numerical parameters of its associated lattice path $\varphi$. For instance if $\varphi$ has domain $I$, then the rank of the underlying poset is $\max I - \min I$, the cohomological degree of the underlying summand is $\max I - \min I - \# I + 1$ and the weight of the underlying summand is $\frac{\max I - \min I - \varphi(\max I) + \varphi(\min I )}{2}$.
\end{myrmk}
\begin{myrmk}\label{rmkklspath}
\phantom\qed
The graded summands of $\Bar(\GGL)$ in $\RKLS$ are exactly the graded summands whose associated lattice paths are strictly positive after $0$, except possibly at the last value. \red{This last value is zero in weight $\rk P / 2$ with $P$ of even rank. Similarly, the graded summands of $\Bar(\GGL)$ in $\iLKLS$ are the graded summands whose associated lattice paths never drop below their last value, except possibly at the first value.}
\end{myrmk}
\begin{nota*}
For any lattice path $\varphi$ \blue{of domain $I$} we denote $\varphi^{\circ} \coloneqq \varphi_{|I \setminus \{\min I, \max I \}}.$
\end{nota*}
Let $M$ be a finite set of integers and $\alpha$ some integer. We denote by $\gr^{M, \alpha} \Bar(\GGL)$ the direct sum of summands of $\Bar(\GGL)$ whose associated lattice path $\varphi$ satisfy
\begin{equation*}
\left\{
\begin{array}{ccc}
    \argmin \varphi^{\circ} &=& M \\
    \min \varphi^{\circ} &=& \alpha,
\end{array}
\right.
\end{equation*} where $\argmin \varphi^{\circ}$ is the set of arguments on which $\varphi^{\circ}$ is minimal (we will call those arguments the internal minima of $\varphi$). Let $S$ be the set of such pairs $(M, \alpha)$. Let us choose any linear order $\leq$ on $S$ satisfying $(M, \alpha) \geq (M', \alpha')$ whenever $\alpha < \alpha'$ or when $\alpha = \alpha '$ and $M \supset M'$. We define an increasing filtration $\filtr^{\bullet}$ indexed by $(S, \leq)$ on $\Bar(\GGL)(P)$ by putting 
\begin{equation*}
    \filtr^{M,\alpha} \Bar (\GGL) \coloneqq \bigoplus_{(M', \alpha') \leq (M, \alpha)} \gr^{M', \alpha'}\Bar(\GGL).
\end{equation*}
One can check that this filtration is compatible with the differential on $\Bar(\GGL)$ (the differential of $\Bar(\GGL)$ sends a summand with associated lattice paths $\varphi$ to summands with associated path $\varphi'$ where $\varphi'$ can be obtained by forgetting one of the arguments of $\varphi$), and that the associated graded of $\filtr^{\bullet}$ is $\gr^{\bullet}$. In the sequel, $\gr^{M, \alpha}$ will mean the associated graded with its differential induced by the differential on $\Bar(\Gerst)$. If $P$ has rank $k$, we have 
\begin{equation*}
    \filtr^{\{1, ... ,k-1\}, 1}\Bar(\GGL)(P) = \RKLS(P).
\end{equation*}
Let $(M = \{x_1< \cdots < x_n\} \neq \emptyset, \alpha \leq 0)$ be an element of $S$. Remark that if $\varphi$ is the lattice path of a summand in $\gr^{M, \alpha}\Bar(\GGL)$ then the lattice path $\varphi_{\llbracket 0, x_{1} \rrbracket}$ only reaches its minimum at $x_1$, or at $x_0$ and $x_1$ (if $\alpha = 0$), \red{and therefore by Remark \ref{rmkklspath} it is the lattice path of some summand in $\iLKLS([G, \un])$} for some element $G$ in $P$ such that $[G, \un]$ has rank $x_1$. Similarly, for any $ 0 \leq i\leq n-1$ the lattice path $\varphi_{\llbracket x_{i}, x_{i+1} \rrbracket}$ only reaches its minimum on both ends and  therefore by Remark \ref{rmkklspath} it is the translation of the lattice path of some summand in some KLS complex of half-weight. Finally, the lattice path $\varphi_{\llbracket x_{n}, k \rrbracket}$ is always strictly above its last value and so it is also the translation of some lattice path of some summand in $\RKLS$. This obviously characterizes all such lattice paths $\varphi$ (see Figure \ref{figpathgraded}). 
\begin{figure}
    \centering
    \begin{tikzpicture}
    \draw [black] (0,0) -- (12,0);
    \node [left=0pt of {(0,0)}, outer sep=0pt] {0};
    \fill [red] (0,0) circle (2.5pt);
    \draw [black] (0,-2) -- (12, -2); 
    \node [left=0pt of {(0,-2)}, outer sep=0pt] {$\alpha$};
    \node [below=0pt of {(3, -2)}, outer sep = 0pt] {$x_1$};
    \node [below=0pt of {(5, -2)}, outer sep = 0pt] {$x_2$};
    \node [below=0pt of {(7, -2)}, outer sep = 0pt] {$x_{n-1}$};
    \node [below=0pt of {(9, -2)}, outer sep = 0pt] {$x_n$};
    \fill [red] (3,-2) circle (2.5pt);
    \fill [red] (5,-2) circle (2.5pt);
    \fill [red] (7,-2) circle (2.5pt);
    \fill [red] (9,-2) circle (2.5pt);
    \fill [red] (12,2) circle (2.5pt);
    
    \draw [red,ultra thick] plot [smooth] coordinates {(0,0) (0.75,0.5) (1.5,-1) (2.25,-0.25) (3,-2)};
    \draw [red,ultra thick] plot [smooth] coordinates {(3,-2) (3.5,-1.5) (4,-1.75) (4.5,-1.25) (5, -2) };
    \draw [red,ultra thick, dashed] plot [smooth, tension = 1.25] coordinates {(5,-2) (6,-1.5) (7,-2)};
    \draw [red,ultra thick] plot [smooth] coordinates {(7,-2) (7.5,-1.5) (8,-1.75) (8.5,-1.25) (9, -2) };
    \draw [red,ultra thick] plot [smooth] coordinates {(9,-2) (10 ,0.25) (11,-1) (12, 2)};
    \draw [black, <->] (0, -2.5) -- (3, -2.5) node[midway,below] {$\iLKLS_{(\star)}$};
    \draw [black, <->] (3, -2.5) -- (5, -2.5) node[midway,below=3pt] {$\RKLS_{\left(\frac{\rk}{2}\right)}$};
    \draw [black, dashed, <->] (5, -2.5) -- (7, -2.5) ;
    \draw [black, <->] (7, -2.5) -- (9, -2.5) node[midway,below=3pt] {$\RKLS_{\left(\frac{\rk}{2}\right)}$};
    \draw [black, <->] (9, -2.5) -- (12, -2.5) node[midway,below=3pt] {$\RKLS_{(\star)}$};
    
\end{tikzpicture}
    \caption{A generic lattice path associated to a summand in $\gr^{M, \alpha}\Bar(\GGL)$.}
    \label{figpathgraded}
\end{figure}

This together with Remark \ref{rmkweightpath} leads to the isomorphism of differential complexes 
\begin{multline}\label{eqgradedisobelow}
     \gr^{M, \alpha}\Bar(\GGL)_{(i)}(P) \cong \\ \bigoplus_{\substack{\zero < G_1 < \cdots < G_n < \un \\ \rk [G_{i-1}, G_{i}] = x_{i} - x_{i-1}}} \left(\RKLS_{\left(\frac{\rk[\zero, G_1] - (k - 2i-\alpha)}{2}\right)}([\zero, G_1])\otimes \RKLS_{\left(\frac{\rk [G_1, G_2]}{2}\right)}([G_1, G_2])\otimes \cdots \right. \\ \left. \otimes \RKLS_{\left(\frac{\rk [G_{n-1}, G_{n}]}{2}\right)}([G_{n-1}, G_{n}]) \otimes
     \iLKLS_{\left(\frac{\rk[G_n, \un]+\alpha}{2}\right)}([G_n, \un])\right).
\end{multline}
(with the convention $x_0 = 0$ and $x_{n+1}=k$). \blue{Note that the intervals in the middle must have even ranks, because the lattice paths we have defined can only go from one value to the same value on a domain of even length.} The identification of the differentials on both sides comes from the fact that the pieces of the differential on some summand of $\gr^{M, \alpha}\Bar(\GGL)(P)$ which forget a point of $M$ land in a strictly smaller filtration grading by definition of $\leq$. Let us compute in which degrees the cohomology of the complex \eqref{eqgradedisobelow} is concentrated. If $n>1$ then the middle KLS complexes in the right hand side of \eqref{eqgradedisobelow} are acyclic by the induction hypothesis and thus the complex \eqref{eqgradedisobelow} is acyclic by the Künneth formula. Likewise, if $\alpha = 0$ then the rightmost KLS complex in the right hand side of \eqref{eqgradedisobelow} is acyclic which implies that the complex \eqref{eqgradedisobelow} is acyclic as well. If $n = 1, \alpha <0$, assume first that we have $i < k/2$. In that case by the induction hypothesis and the Künneth formula the above complex is concentrated in degree 
\begin{align*}
    \frac{x_1 + \alpha}{2}+ \frac{k-x_1 - (k - 2i-\alpha)}{2} &= i + \alpha \\
    &\leq i-1. 
\end{align*} 
Let us assume now that we have $i> k/2$. If $\alpha = k-2i$ then the leftmost KLS complex in the right hand side of \eqref{eqgradedisobelow} is acyclic which means that the complex \eqref{eqgradedisobelow} is acyclic. Otherwise if $k-2i > \alpha $ then the cohomology of the complex \eqref{eqgradedisobelow} is concentrated in degree 
\begin{align*}
    \frac{x_1 + \alpha}{2}+ \frac{k-x_1 - (k - 2i-\alpha)}{2} &= i + \alpha \\
    &\leq i-2
\end{align*} 
($\alpha \leq -2$ because we have $\alpha < k-2i < 0$). Finally if $k-2i < \alpha$ then the cohomology of the complex \eqref{eqgradedisobelow} is concentrated in degree 
\begin{align*}
    \frac{x_1 + \alpha}{2}+ \frac{k-x_1 - (k - 2i-\alpha)}{2} - 1&= i + \alpha -1 \\
    &\leq i-2.
\end{align*} 

On the other hand if $(M = \{x_1< \cdots < x_n\} \neq \emptyset, \alpha)$ is an element of $S$ such that $\alpha$ is strictly positive, then by the same arguments as above we have the same isomorphism of differential complexes \eqref{eqgradedisobelow}. As previously, by our induction hypothesis and the Künneth formula one can compute in which degrees the cohomology of the complex \eqref{eqgradedisobelow} is concentrated. If $i < k/2$ and $\alpha < k-2i$ we get cohomological degree 
\begin{align*}
    \frac{x_1+\alpha}{2} - 1 + \frac{k-x_1-(k-2i-\alpha)}{2} &=  i + \alpha - 1 \\
    &\geq i. 
\end{align*}
If $i < k/2$ and $\alpha > k-2i$ we get cohomological degree 
\begin{align*}
    \frac{x_1+\alpha}{2} - 1 + \frac{k-x_1-(k-2i-\alpha)}{2} - 1  &=  i + \alpha - 2 \\
    &\geq i 
\end{align*}
($\alpha$ is greater or equal than 2 because we have $\alpha > k-2i > 0$). Otherwise, if $i > k/2$ we get cohomological degree 
\begin{align*}
    \frac{x_1+\alpha}{2} - 1 + \frac{k-x_1-(k-2i-\alpha)}{2} - 1 &=  i + \alpha - 2 \\
    &\geq i-1. 
\end{align*}
The case of $\LKLS$, $\iLKLS$ and $\iRKLS$ being completely symmetric, this concludes the induction step of Statement i). 

One can consider a coarsening $\F'$ of the filtration $\F$ obtained by only looking at the height $\alpha$ of the internal minima (and not at the minima themselves). Let us assume that we are in the case $k-2i = 1$ (keeping the same notations as in the proof of Statement i). By the computations carried out previously the first page of the spectral sequence associated to $\F'$ restricted to $\RKLS$ has the shape depicted in Figure \ref{figE1RKLS}. In this figure the vertical axis is given by the filtration grading (the height of the internal minima) and the horizontal axis is given by the cohomological degree. 
\begin{figure}[H]
\center
\begin{tikzcd}
  0 &  0 & 0 & ... \\
  0 &  0 & A_3 \ar[ur, "d_1"] & 0\\
  0 &  A_2 \ar[ur, "d_1"] & 0 & 0\\
  0 & 0 & 0 & 0 
\end{tikzcd}
\caption{The first page of the spectral sequence associated to $\F'_{|\RKLS}$ in the case $k-2i = 1$.}
\label{figE1RKLS}
\end{figure}
This spectral sequence stabilizes at the second page and the cohomology of $\RKLS_{(i)}$ is isomorphic to the kernel of the leftmost differential $d_1$. The leftmost term $A_2$, living in filtration grading $2$, is given by 
\begin{equation*}
A_2 =   \bigoplus_{\substack{j < k \\ G \in P \textrm{ s.t. } \rk G = j}} \left( H^{\frac{k-j+1}{2}-1}\left(\RKLS_{\left(\frac{k-j+1}{2}\right)}([\zero, G])\right) \otimes H^{\frac{j+2}{2} - 1}\left( \iLKLS_{\left(\frac{j+2}{2}\right)}([G, \un])\right)  \right).   
\end{equation*}

If on the other hand we have $k-2i \geq 2$ then the first page looks slightly different because the first non trivial term lives in filtration grading $1$ and there will be a jump across filtration grading $k-2i$ (the graded complex is acyclic in this filtration grading). In this case the leftmost term $A_1$ is given by 
\begin{equation*}
A_1 =   \bigoplus_{\substack{j < k \\ G \in P \textrm{ s.t. } \rk G = j}} \left(H^{\frac{2i - j +1}{2}}\left(\RKLS_{\left(\frac{2i-j+1}{2}\right)}([\zero, G])\right) \otimes H^{\frac{j+1}{2} - 1}\left( \iLKLS_{\left(\frac{j+1}{2}\right)}([G, \un])\right)  \right).   
\end{equation*}
The spectral sequence only stabilizes at the third page but the cohomology of $\RKLS_{(i)}$ can still be identified with the kernel of some differential defined on $A_1$. Putting those two cases together we get the following lemma which we will need in the proof of Statement ii).
\begin{mylemma}\label{lemmarepresentation}
    If $k-2i = 1$ (resp. $k-2i \geq 2$) then any element in $H^i(\RKLS_{(i)})$ can be represented by a sum of homogeneous elements in graded summands whose associated lattice path have only one internal minimum of height $2$ (resp. height $1$), and similarly for $\iLKLS$.  
\end{mylemma}
For Statement ii) let us assume that the rank $k$ is even. We depict the first page of the spectral sequence associated to the coarsening $\F'$ on the whole Bar complex of weight $k/2$ in Figure \ref{figE1KLShalf}. As in Figure \ref{figE1RKLS} the vertical axis is the filtration grading and the horizontal axis is the cohomological degree. 

\begin{figure}[H]
\center
\begin{tikzcd}
   &  & ... &  & ... \\
   & 0 & 0 & A_2 \ar[ur, "d_1"] & \\
   & 0 & A_1 \ar[ur, "d_1"] & 0 & \\
  ... & 0 & 0 & 0 & ...\\
   & 0 & A_{-1} \ar[uuur, dashed, "d_3"'] & 0 & \\
   & A_{-2} \ar[uuur, dashed, "d_3"] \ar[ur, "d_1"'] & 0&0 & \\
  ... \ar[ur, "d_1"'] & &... & & 
\end{tikzcd}
\caption{The first page of the spectral sequence associated to $\F'$ in half weight.}
\label{figE1KLShalf}
\end{figure}

The term $A_{-1}$ is given by  
\begin{equation*}
    A_{-1} = \bigoplus_{\substack{j < k \\ G \in P \textrm{ s.t. } \rk G = j}} \left(H^{\frac{k-j-1}{2}}\left(\RKLS_{\left(\frac{k-j-1}{2}\right)}([\zero, G])\right) \otimes H^{\frac{j-1}{2}}\left( \iLKLS_{\left(\frac{j-1}{2}\right)}([G, \un])\right) \right).
\end{equation*}  
By Lemma \ref{lemmarepresentation} an element of $A_{-1}$ can be represented as a sum of homogeneous elements whose associated lattice path is the concatenation of two lattice paths each having only one internal minimum of height $1$, and the concatenation itself having only one internal minimum of height $-1$. The differential of such homogeneous elements has no component landing in a graded summand whose associated lattice path has only one internal minimum of height $2$, which proves that the differential $d_3$ is equal to zero on $A_{-1}$. By a similar argument it is also zero on $A_{-2}$. This implies that the spectral sequence has already stabilized at the second page, to zero above $A_1$ by the acyclicity of $\Bar(\GGL)$ in strictly positive degrees. However, the second page of the spectral sequence above $A_1$ is also the second page of the spectral sequence associated to the restriction of the filtration $\F$ to $\RKLS_{(\frac{k}{2})}$. This proves the acyclicity of $\RKLS_{(\frac{k}{2})}$. The case of $\LKLS$, $\iRKLS$ and $\iLKLS$ being completely symmetric, this concludes the induction step of Statement ii) and the proof of Theorem \ref{theomainhomo}.

\end{proof}
We end this article with a series of informal discussions aimed at further research. 
\begin{discussion}[What about equivariant Kazhdan--Lusztig--Stanley theory ?]
The operadic framework developed in this article should also be suitable to handle equivariant Kazhdan--Lusztig--Stanley theory (see \cite{braden2023singular} Appendix A for a reference on this topic). One just needs to add the group actions as part of the datum of the $\P$-operads, which was originally considered in \cite{coron2023matroids} for geometric lattices. The cooperad $(\GGL^{\GL})^{\vee}$ has an equivariant enhancement with automorphism group action defined by 
\begin{equation*}
\begin{array}{ccc}
  \GGL(\L)^{\vee} \simeq \OS(\L) & \xrightarrow{\GGL^{\vee}(\phi)} & OS(\L')\simeq \GGL(\L')^{\vee} \\
        e_H & \rightarrow & e_{\phi^{-1}(H)}
\end{array}
\end{equation*}
for any isomorphism $\phi: \L' \rightarrow \L$ between geometric lattices. This leads to automorphism group action on bar construction and ultimately on KLS complexes. 
\end{discussion}
\begin{discussion}[What about Hodge theory ?]
It would be interesting to relate the Hodge theoretic methods of \cite{braden2023singular} to the methods of  this article. Many of the protagonists appearing in \cite{braden2023singular}, such as augmented and non augmented Chow rings of geometric lattices, or Rouquier complexes, have an operadic interpretation (the non augmented Chow rings have a $\GL$-cooperadic structure studied in \cite{coron2023matroids}, the augmented Chow rings have a structure of an operadic comodule over the non augmented Chow rings, and the Rouquier complexes can be interpreted as bar constructions of those operadic structures). At the moment the material in this article cannot account for any Hodge theoretic result about KLS polynomials, because we have no structure that relates KLS complexes of different weight. This could be remedied by considering $\GGL$ as an operad in coalgebras instead of just vector spaces. 
\end{discussion}
\begin{discussion}[What about geometry ?]
A natural question to ask is whether there is a geometric operadic structure behind the algebraic operadic structure $\Gerst$, i.e. an operad in geometric objects whose homology is given by $\Gerst$, when restricting to realizable geometric lattices. As mentioned in the introduction, for braid arrangements the answer is given by the little 2-discs operad, but this generalizes very poorly to other hyperplane arrangements over $\Cx$. Fortunately, there exists a multitude of other geometric operads which give the operad $\Gerst$ after passing to homology (such operads are called $E_2$-operads). One of those $E_2$ operads is given by the real Fulton--MacPherson compactifications of braid arrangements (see \cite{getzler1994operads}). We conjecture that this new candidate generalizes to any hyperplane arrangement over $\Cx$ to give a ``geometrification'' of the $\GL$-operad $\Gerst$ (restricted to geometric lattices realizable over $\Cx$).
\end{discussion}

\textbf{Acknowledgements.} The author would like to thank Alex Fink for his availability and many useful conversations. Many thanks as well to the tropical and geometric combinatorics team members for their warm welcome at Queen Mary University of London. This work was supported by the Engineering and Physical Sciences Research Council [grant number EP/X001229/1].

\bibliographystyle{apalike}
\bibliography{sample}

\end{document}